\documentclass[final,1p,times]{elsarticle}
\usepackage{amsmath, amsfonts, amssymb}
\usepackage{graphicx,epsfig}\graphicspath{{images/}}
\usepackage[neverdecrease]{paralist}
\numberwithin{equation}{section}

\usepackage{dsfont}
\usepackage{float}
\usepackage{array}
\usepackage{multirow}
\usepackage{booktabs}
\usepackage{color,colortbl}

\usepackage{amsthm}
\theoremstyle{plain}\newtheorem{remark}{Remark}[section]
\theoremstyle{definition}

\DeclareMathOperator{\argmax}{argmax}

\DeclareMathOperator{\std}{std}

\definecolor{shadegray}{rgb}{0.8,0.8,0.8}

\journal{Applied and Computational Harmonic Analysis}

\begin{document}

\begin{frontmatter}

\title{Nonlinear Mode Decomposition: \\a new noise-robust, adaptive decomposition method\tnoteref{support}}
\tnotetext[support]{Work supported by the Engineering and Physical Sciences Research Council, UK}

\author{Dmytro Iatsenko}
\ead{dmytro.iatsenko@gmail.com}

\author{Peter V. E. McClintock\corref{}}
\ead{p.v.e.mcclintock@lancaster.ac.uk}

\author{Aneta Stefanovska\corref{cor1}}
\cortext[cor1]{Corresponding author}
\ead{aneta@lancaster.ac.uk}

\address{Department of Physics, Lancaster University, Lancaster LA1 4YB, UK}

\begin{abstract}
We introduce a new adaptive decomposition tool, which we refer to as {\it Nonlinear Mode Decomposition} (NMD). It decomposes a given signal into a set of physically meaningful oscillations for any waveform, simultaneously removing the noise. NMD is based on the powerful combination of time-frequency analysis techniques -- which together with the adaptive choice of their parameters make it extremely noise-robust -- and surrogate data tests, used to identify interdependent oscillations and to distinguish deterministic from random activity. We illustrate the application of NMD to both simulated and real signals, and demonstrate its qualitative and quantitative superiority over the other existing approaches, such as (ensemble) empirical mode decomposition, Karhunen-Lo\`{e}ve expansion and independent component analysis. We point out that NMD is likely to be applicable and useful in many different areas of research, such as geophysics, finance, and the life sciences. The necessary MATLAB codes for running NMD are freely available for download.
\end{abstract}


\end{frontmatter}

\section{Introduction}\label{sec:introduction}

In real life, measured signals are often composed of many different oscillations (or modes), with complex waveforms and time-varying amplitudes and frequencies. These oscillations contain a great deal of valuable information about the originating system and they are therefore worthy of very careful study. Their properties can be used e.g.\ to predict earthquakes from geophysical signals \cite{Beroza:90}, or to assess health from signals generated by the human body \cite{Iatsenko:13a,Shiogai:10}. However, to gain access to all of the properties of a particular mode, it should be first extracted from the signal. For a reliable analysis, therefore, one should decompose the signal, i.e.\ recover the individual oscillations that are present in it, separating them both from each other and from the inevitable background of noise.

How best to accomplish this decomposition is a problem of long-standing, and many different approaches have been proposed. However, to the best of the authors knowledge, each of the current methods has at least two of the following flaws:
\begin{enumerate}
\item The method contains user-defined parameters and is quite sensitive to their choice. These parameters cannot be adaptively chosen and there are no more-or-less universal settings.
\item The method is not noise-robust.
\item If an individual mode has a complex (non-sinusoidal) waveform, the method will decompose it into a few oscillations with simpler waveforms.
\item The modes returned by the method are not always physically meaningful. For example, the method will decompose even a random signal, such as Brownian noise, onto a set of oscillations.
\end{enumerate}

\noindent The popular Empirical Mode Decomposition (EMD) \cite{Huang:98} method suffers from disadvantages 2, 4 and 3, with the latter manifesting itself only in cases where the corresponding oscillation has a waveform with multiple peaks. EMD is very sensitive to noise and, to overcome this limitation, Wu and Huang \cite{Wu:09} recently proposed a variant called ensemble empirical mode decomposition (EEMD). The idea is to add independent realizations of white Gaussian noise to the signal and apply EMD each time, with the final modes being obtained as the ensemble-averages of the modes for each noise realization. Even so, however, EEMD still cannot be regarded as noise-robust method (at least in comparison with other methods), as will be seen in the examples below. Additionally, it contains a non-adaptive parameter (drawback 1) controlling the standard deviation of the added noise, though its choice does not seem to exert a strong influence on the results \cite{Wu:09}.

For multivariate time-series the most widely-used decomposition methods are principal component analysis (PCA) \cite{Jolliffe:02,Abdi:10} and independent component analysis (ICA) \cite{Comon:94,Hyvarinen:00}. In the present work, however, we restrict our consideration to univariate signals. To apply PCA and ICA in this case, one should first construct a multivariate signal from the given single time-series, which can be accomplished by using its time-shifted blocks. PCA and ICA can then be applied to this collection of original signal parts, in which case they are called single-channel PCA (also known as Karhunen-Lo\`{e}ve expansion) and single-channel ICA \cite{Davies:07}, respectively. However, both of these methods are generally susceptible to the time-shifts and the number of blocks used to create the multivariate signal from the original one (as was demonstrated e.g.\ in \cite{Hozic:00} for the Karhunen-Lo\`{e}ve expansion), thus suffering from the drawback 1; in addition, they also contain flaw number 4.

Another approach, which is now becoming increasingly popular, is to use a particular time-frequency representation (TFR) \cite{Iatsenko:13tfr1,Mallat:08,Addison:10,Cohen:95,Hlawatsch:92,Daubechies:92} of the signal, e.g.\ the wavelet transform, for decomposing it into separate modes. Qualitatively, the signal is first mapped onto the time-frequency plane, where its oscillatory components appear in the form of ``curves'', formed by the sequences of amplitude peaks. Having identified these curves, one can then recover the corresponding oscillations, thereby decomposing the signal. Many different techniques for realizing this idea have been proposed \cite{Delprat:92,Iatsenko:13ridge,Carmona:99,Carmona:97,Daubechies:11,Thakur:13}, with the great advantage of such approaches being their high noise-robustness \cite{Iatsenko:13tfr2,Chen:12}. However, apart from the disadvantage 4, which is common to all current methods, any TFR-based decomposition also suffers from drawback 3, because any oscillation with a non-sinusoidal waveform is typically represented by a number of curves in the time-frequency plane.

The flaws mentioned greatly restrict the applicability of the approaches currently in use, so that only for a very small class of signals can the decomposition be carried out successfully. To overcome these limitations, we now introduce a new method -- Nonlinear Mode Decomposition (NMD) -- which is free from all the drawbacks considered above. It is based on the powerful combination of the time-frequency analysis techniques reviewed/developed in \cite{Iatsenko:13tfr1,Iatsenko:13tfr2,Iatsenko:13ridge}, surrogate data tests \cite{Schreiber:00,Theiler:92,Andrzejak:03,Quiroga:02}, and the idea of harmonic identification \cite{Sheppard:11}. NMD has proven to be extremely noise-robust, and it returns only the physically meaningful oscillations with any waveforms, at the same time removing the noise. Furthermore, we develop a set of criteria using which almost all NMD settings can be adapted automatically to the signal, greatly improving its performance and making it a kind of super-adaptive approach. These features make NMD unique, providing many significant advantages over other existing approaches, as we demonstrate below on both simulated and real data.

The structure of the work is as follows. In Sec.\ \ref{sec:background} we summarize the notation used and review the necessary background. The basic idea of the NMD and its implementation are considered in Sec.\ \ref{sec:nmd}. In Sec.\ \ref{sec:improve} we introduce some important upgrades to improve the NMD performance, while the choice of its parameters is considered in Sec.\ \ref{sec:pchoice}. We demonstrate the power of the method by applying it to simulated signals in Sec.\ \ref{sec:simex}, and to real signals in Sec.\ \ref{sec:realex}. We conclude in Sec.\ \ref{sec:conclusions}.

\section{Background and terminology}\label{sec:background}

\subsection{Main notation}

Given some function $f(t)$, we will denote its Fourier transform, positive frequency part, time-average and standard deviation as $\hat{f}(\xi)$, $f^{+}(t)$, $\langle f(t)\rangle$ and $\std[f(t)]$, respectively:
\begin{equation}\label{nt}
\begin{gathered}
\hat{f}(\xi)\equiv\int_{-\infty}^{\infty} f(t)e^{-i\xi t}dt\Leftrightarrow
f(t)=\frac{1}{2\pi}\int_{-\infty}^{\infty}\hat{f}(\xi)e^{i\xi t}d\xi,\\
f^{+}(t)\equiv \int_{0^+}^\infty \hat{f}(\xi)e^{i\xi t}d\xi,\;
\langle f(t)\rangle=\frac{1}{T}\int_0^T f(t)dt,\;
\std[f(t)]\equiv\sqrt{\langle [f(t)]^2\rangle-[\langle f(t)\rangle]^2}\\
\end{gathered}
\end{equation}
where $T$ denotes the overall time-duration of $f(t)$ (in theory $T\rightarrow\infty$).

Next, the \emph{AM/FM component} (or simply \emph{component}) is defined as the function of time $t$ having the form
\begin{equation}\label{amfm}
x(t)=A(t)\cos\phi(t),\;(A(t)>0,\;\nu(t)=\phi'(t)>0\;\forall t),
\end{equation}
where $A(t)$, $\phi(t)$ and $\nu(t)$ are respectively the instantaneous amplitude, phase and frequency of the component. We restrict the consideration to components satisfying the analytic approximation
\begin{equation}\label{aa}
[A(t)\cos\phi(t)]^+\approx \frac{1}{2}A(t)e^{i\phi(t)}
\end{equation}
because, otherwise, the amplitude and phase are not uniquely defined and cannot in principle be recovered reliably by time-frequency analysis methods \cite{Iatsenko:13tfr1}. For a more detailed discussion of the definitions of $A(t),\phi(t),\nu(t)$ and of some related issues see e.g.\ \cite{Picinbono:97,Boashash:92a,Boashash:92b}.

In real cases, a signal usually contains many components of the form (\ref{amfm}). Moreover, some of them typically do not correspond to an independent activity, but arise due to complex waveform of a particular mode to which they are related (see below). This is because real oscillations are rarely purely sinusoidal, but have more complicated shapes as the result of nonlinearities in the generating system and/or the measurement apparatus. For example, the AM/FM component (\ref{amfm}) raised to the third power $[A(t)\cos\phi(t)]^3=\frac{3}{4}A^3(t)[\cos\phi(t)+\frac{1}{3}\cos 3\phi(t)]$ consists already of two components, although there is obviously only one meaningful oscillation.

It is therefore better to consider the full \emph{Nonlinear Modes}, defined as the sum of all components corresponding to the same activity:
\begin{equation}\label{nm}
c(t)=A(t){\rm v}(\phi(t))=A(t)\sum_h a_h\cos(h\phi(t)+\varphi_h).
\end{equation}
where ${\rm v}(\phi(t))={\rm v}(\phi(t)+2\pi)$ is some periodic function of phase (also called the ``wave-shape function'' \cite{Wu:13}), which due to its periodicity can always be expanded in the Fourier series (\ref{nm}). Without loss of generality, we fix the normalization of $A(t)$ and $\phi(t)$ in (\ref{nm}) by setting $a_1=1$, $\varphi_1=0$. Then the instantaneous phase and frequency of the whole NM are $\phi(t)$ and $\nu(t)=\phi'(t)$, respectively \cite{Wu:13}. In what follows, the AM/FM components composing the NM will be referred to as \emph{harmonics}, with the $h^{\rm th}$ harmonic being represented by a term $\sim \cos(h\phi(t)+\varphi_h)$ in (\ref{nm}).

We assume that the signal is composed of the NMs $c_i(t)$ of form (\ref{nm}) plus some noise $\eta(t)$ (the class of noise will be considered later, in Sec.\ \ref{sec:stopcrit}):
\begin{equation}\label{sig}
s(t)=\sum_i c_i(t)+\eta(t).
\end{equation}
Our ultimate goal is then to extract all the NMs present in the signal and to find their characteristics, such as their amplitudes $A(t)$, phases $\phi(t)$ and frequencies $\nu(t)$, as well as the amplitude scaling factors $a_h$ and phase shifts $\varphi_h$ of the harmonics.

\subsection{Time-frequency representations}\label{sec:TimeFrequencyRepresentations}

The time-frequency representation (TFR) of the signal is constructed by projecting it onto the time-frequency plane, which allows us to visualize and follow the evolution of the signal's spectral content in time. Depending on how the projection is constructed, there are many different types of TFRs. In our case, however, we will need those from which the components present in the signal can be extracted and reconstructed in a straightforward way. TFRs meeting this criterion are the windowed Fourier transform (WFT) $G_s(\omega,t)$ and the wavelet transform (WT) $W_s(\omega,t )$. All aspects of their use in practice, and our implementation of them, are discussed in \cite{Iatsenko:13tfr1} (see also the classical reviews and books \cite{Mallat:08,Addison:10,Cohen:95,Hlawatsch:92,Daubechies:92}), to which we refer the reader for detail; here we just review briefly the minimal necessary background.

Given a real signal $s(t)$, its WFT and WT are constructed as
\begin{equation}\label{wftwt}
\begin{aligned}
&G_s(\omega,t)\equiv\int_{-\infty}^\infty s^{+}(u)g(u-t)e^{-i\omega (u-t)}du
=\frac{1}{2\pi}\int_0^\infty e^{i\xi t}\hat{s}(\xi)\hat{g}(\omega-\xi)d\xi,\\
&W_s(\omega,t)\equiv\int_{-\infty}^\infty s^{+}(u)\psi^*{\Big(}\frac{\omega(u-t)}{\omega_\psi}{\Big)}\frac{\omega du}{\omega_\psi}
=\frac{1}{2\pi}\int_0^\infty e^{i\xi t}\hat{s}(\xi)\hat{\psi}^*\Big(\frac{\omega_{\psi}\xi}{\omega}\Big)d\xi,
\end{aligned}
\end{equation}
where $g(t)$ and $\psi(t)$ are respectively the chosen window and wavelet functions, with the latter satisfying the admissibility condition $\hat{\psi}(0)=0$; the $\omega_\psi=\operatorname{argmax}|\hat{\psi}(\xi)|$ denotes wavelet peak frequency, while for the WFT we assume $\hat{g}(\xi)$ to be centered around its peak: $|\hat{g}(0)|=\max|\hat{g}(\xi)|$. We choose a Gaussian window for the WFT and a lognormal wavelet for the WT:
\begin{equation}\label{winwav}
\begin{aligned}
&\hat{g}(\xi)=e^{-(f_0\xi)^2/2}\Leftrightarrow g(t)=\frac{1}{\sqrt{2\pi}f_0}e^{-(t/f_0)^2/2}\\
&\hat{\psi}(\xi)=e^{-(2\pi f_0\log\xi)^2/2},\;\omega_\psi=1\\
\end{aligned}
\end{equation}
where $f_0$ is the resolution parameter determining the tradeoff between time and frequency resolution of the resultant transform (by default we are using $f_0=1$). The forms (\ref{winwav}) appear to be the most useful ones in terms of the joint time-frequency resolution and the reconstruction of components \cite{Iatsenko:13tfr1}. Note, however, that all of the considerations that follow are applicable for any $\hat{g}(\xi)$ and $\hat{\psi}(\xi>0)$ which are real, positive and unimodal around the single peak.

The difference between the WFT and WT lies in the way they resolve the components present in the signal. The WFT has linear frequency resolution, i.e.\ distinguishes components on the basis of their frequency differences, whereas the WT has logarithmic frequency resolution, distinguishing components on the basis of ratios (or the differences between the logarithms) of their frequencies. As a result, the time resolution of the WT increases with frequency, so that the time variability of the higher frequency components is represented better than that for the components of lower frequency, while the time resolution of the WFT is the same at all frequencies.

In theory, the WFT and WT (\ref{wftwt}) depend on continuous $\omega$, but in practice they are calculated at particular discrete $\omega_k$. The latter can be chosen as $\omega_k=(k-k_0)\Delta\omega$ for the WFT and $\log\omega_k=(k-k_0)\frac{\log 2}{n_v}$ for the WT, reflecting the linear and logarithmic frequency resolution of the corresponding TFRs; the discretization parameters $\Delta\omega$ and $n_v$ can be selected using the criteria established in \cite{Iatsenko:13tfr1}. Due to the finite time-length of real signals, all TFRs suffer from boundary distortions and, to reduce them, we use the predictive padding suggested in \cite{Iatsenko:13tfr1}.

\section{Nonlinear Mode Decomposition: basic idea and ingredients}\label{sec:nmd}

The main goal of NMD is to decompose a given signal into a set of Nonlinear Modes (\ref{nm}). To do this, four steps are necessary:

\begin{itemize}

\item[(a)] Extract the fundamental harmonic of an NM accurately from the signal's TFR.

\item[(b)] Find candidates for all its possible harmonics, based on its properties.

\item[(c)] Identify the true harmonics (i.e.\ corresponding to the same NM) among them.

\item[(d)] Reconstruct the full NM by summing together all the true harmonics; subtract it from the signal, and iterate the procedure on the residual until a preset stopping criterion is met.

\end{itemize}

\noindent These individual sub-procedures are explained in detail in the Sections below.

Note that NMD can be based either on the WFT, or on the WT, although in what follows (Sec.\ \ref{sec:ptfrtype}) we will combine them to make the procedure more adaptive. Nevertheless, for the time being we assume that we have selected one of these TFR types and that we will use it for all operations.

\subsection{Component extraction and reconstruction}\label{sec:cextract}

\begin{figure}[t]
\includegraphics[width=1.0\linewidth]{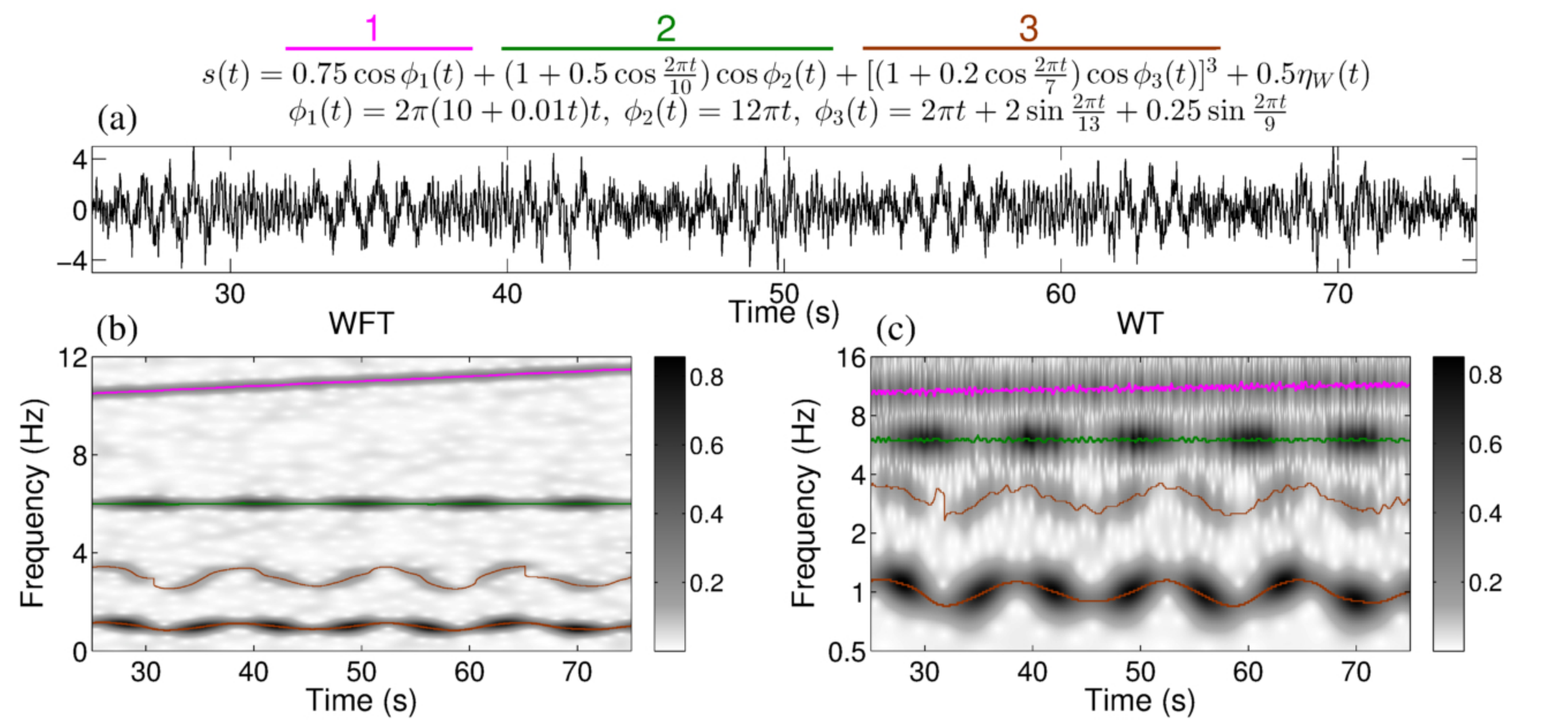}
\caption{(a) The central 50\,s section of the signal $s(t)$ specified at the top, which is composed of three oscillations corrupted with white Gaussian noise of standard deviation equal to $0.5$. The oscillations are: 1) near 10\,Hz, with frequency modulation (chirp); 2) near 6\,Hz, with amplitude modulation; 3) near 1\,Hz, with both modulations and a non-sinusoidal waveform. (b),(c) Respectively, the central parts of the WFT and WT of the signal $s(t)$ shown in (a), where the ridge curves related to each oscillation are shown by solid lines of the corresponding colors; the third oscillation is represented by two curves because of its complex waveform, while the apparent distortions in the ridge frequency profiles are due to noise. The signal was sampled at 100\,Hz for 100\,s.}
\label{fig:ex0tfr}
\end{figure}

Provided that the TFR has sufficient resolution in time and frequency to resolve components that are close in frequency and, at the same time, to portray reliably their amplitude/frequency modulation, each component will be represented by a unique sequence of TFR amplitude peaks (also called ridge points). We will call such a sequence a \emph{ridge curve}, and denote it as $\omega_p(t)$; an illustrative example showing the WFT and WT of a signal and the ridge curves of its components is presented in Fig.\ \ref{fig:ex0tfr}. Note that oscillations with a non-sinusoidal waveform will appear as two or more components in a TFR, consisting of the fundamental and its harmonics.

To extract particular a component from a TFR, one needs first to find its ridge curve. This is a non-trivial task, because in real cases it is not always clear which sequence of peaks corresponds to which component, and which peaks are just noise-induced artifacts. The procedures for ridge curve extraction are discussed and developed in \cite{Iatsenko:13ridge}. By default we use here scheme II(1,1) from the latter, as it appears to be very universal in applicability. In short, in the case of WFT we select the ridge curve $\omega_p(t)$ as being that sequence of peaks (among all those that are possible) which maximizes the path functional
\begin{equation}\label{ecurve}
\begin{aligned}
S[\omega_p(t)]\equiv \sum_{n=1}^N\bigg[\log|G_s(\omega_p(t_n),t_n)|
-\frac{1}{2}\bigg(\frac{|\Delta\omega_p(t_n)-\langle\Delta\omega_p\rangle|}{\std[\Delta\omega_p]}\bigg)^2
-\frac{1}{2}\bigg(\frac{\omega_p(t_n)-\langle\omega_p\rangle}{\std[\omega_p]}\bigg)^2\bigg]
\end{aligned}
\end{equation}
where $\Delta\omega_p(t_n)\equiv\omega_p(t_n)-\omega_p(t_{n-1})$, while $\langle f(t)\rangle\equiv (1/N)\sum f(t_n)$ and $\std[f(t)]\equiv[\langle (f(t))^2\rangle-\langle f(t)\rangle^2]^{1/2}$ denote the time-average and standard deviation, respectively. In the case of the WT, one uses the same functional (\ref{ecurve}), except that now all is taken on a logarithmic frequency scale ($\omega_p\rightarrow\log\omega_p$, $\Delta\omega_p\rightarrow\Delta\log\omega_p$). The approximate optimal $\omega_p(t)$ in terms of (\ref{ecurve}) can be found in $O(N)$ steps as described in \cite{Iatsenko:13ridge}.

Having found the ridge curve $\omega_p(t)$, the amplitude $A(t)$, phase $\phi(t)$ and frequency $\nu(t)\equiv\phi'(t)$ of the corresponding component $x(t)=A(t)\cos\phi(t)$ can be reconstructed by either of two methods, ridge and direct \cite{Iatsenko:13tfr1}. The ridge method utilizes only the TFR values at the peak, and the corresponding formulas read
\begin{equation}\label{ridgerec}
\begin{aligned}
\mbox{\textbf{ridge[WFT]:}}&\;
\nu(t)=\omega_p(t)+\delta\nu_d(t),\;
A(t)e^{i\phi(t)}=\frac{2G_s(\omega_p(t),t)}{\hat{g}(\omega_p(t)-\nu(t))},\\
\midrule
\mbox{\textbf{ridge[WT]:}}&\;
\nu(t)=\omega_p(t)e^{\delta\log\nu_d(t)},\;
A(t)e^{i\phi(t)}=\frac{2W_s(\omega_p(t),t)}{\hat{\psi}^*(\omega_\psi\nu(t)/\omega_p(t))},
\end{aligned}
\end{equation}
where $\delta \nu_d(t)$ and $\delta\log\nu_d(t)$ are the corrections for discretization effects. Thus, in theory one has a continuous frequency scale, so that the peak positions $\omega_p(t)$ are exact. In practice, however, TFRs are calculated at the discrete values $\omega_k$, thereby introducing errors into the parameter estimates. To reduce these errors, one can use quadratic interpolation to better locate the peak positions (that are also used subsequently to improve the accuracy of the amplitude estimates in (\ref{ridgerec})), with the corresponding correction terms being
\begin{equation}\label{deltanu}
\begin{aligned}
\mbox{\textbf{WFT:}}&\;\delta\nu_d(t)=\frac{\Delta\omega}{2}\frac{a_3-a_1}{2a_2-a_1-a_3},\;
a_{\{1,2,3\}}\equiv |G_s(\omega_{\{k_p(t)-1,k_p(t),k_p(t)+1\}},t)|,\\
\mbox{\textbf{WT:}}&\;\delta\log\nu_d(t)=\frac{n_v^{-1}\log 2}{2}\frac{a_3-a_1}{2a_2-a_1-a_3},\;
a_{\{1,2,3\}}\equiv |W_s(\omega_{\{k_p(t)-1,k_p(t),k_p(t)+1\}},t)|,\\
\end{aligned}
\end{equation}
where $k_p(t)$ denote the discrete index of the peak at each time: $\omega_p(t)\equiv\omega_{k_p(t)}$.

Another way is to reconstruct the component from the full time-frequency region into which it is mapped for a given TFR. We call this region the \emph{time-frequency support} and denote it as $[\omega_-(t),\omega_+(t)]$. It can be defined as the widest region of unimodal and non-zero TFR amplitude around the ridge point $\omega_p(t)$ at each time \cite{Iatsenko:13tfr1}. Thus, at $\omega_{\pm}(t)$ the TFR amplitude either becomes zero or starts growing with increasing distance from the peak $|\omega-\omega_p(t)|$. The parameters of the component can be reconstructed from its time-frequency support as
\begin{equation}\label{dirrec}
\begin{aligned}
&\mbox{\textbf{direct[WFT]:}}\\
&A(t)e^{i\phi(t)}=C_g^{-1}\int_{\omega_-(t)}^{\omega_+(t)} G_s(\omega,t)d\omega,\;
C_g\equiv\frac{1}{2}\int_{-\infty}^{\infty} \hat{g}(\xi)d\xi=\pi g(0)\\
&\nu(t)=\operatorname{Re}\Bigg[\frac{\int_{\omega_-(t)}^{\omega_+(t)} \omega G_s(\omega,t)d\omega}
{\int_{\omega_-(t)}^{\omega_+(t)} G_s(\omega,t)d\omega}-\overline{\omega}_g\Bigg],\;
\overline{\omega}_g\equiv\frac{C_g^{-1}}{2}\int_{-\infty}^{\infty} \xi\hat{g}(\xi)d\xi,\\
&(\overline{\omega}_g=0\mbox{ for symmetric }\hat{g}(\xi),\mbox{ e.g.\ Gaussian window}),\\\midrule
&\mbox{\textbf{direct[WT]:}}\\
&A(t)e^{i\phi(t)}=C_\psi^{-1}\int_{\omega_-(t)}^{\omega_+(t)} W_s(\omega,t)\frac{d\omega}{\omega},\;
C_\psi\equiv\frac{1}{2}\int_0^\infty\hat{\psi}^*(\xi)\frac{d\xi}{\xi},\\
&\nu(t)=\operatorname{Re}\frac{D_\psi^{-1}\int_{\omega_-(t)}^{\omega_+(t)} \omega W_s(\omega,t)\frac{d\omega}{\omega}}
{C_\psi^{-1}\int_{\omega_-(t)}^{\omega_+(t)} W_s(\omega,t)\frac{d\omega}{\omega}},\;
D_\psi\equiv\frac{\omega_\psi}{2}\int_0^\infty \frac{1}{\xi}\hat{\psi}^*(\xi)\frac{d\xi}{\xi}.\\
\end{aligned}
\end{equation}

\noindent The difference and relative performance of the two reconstruction methods were studied in \cite{Iatsenko:13tfr2}. When the component is perfectly represented in the TFR, i.e.\ $[\omega_-(t),\omega_+(t)]$ contain all of its power and no side influences, direct reconstruction is by definition exact (up to the error of the analytic signal approximation, see \cite{Iatsenko:13tfr1}), while the ridge estimates contain an error proportional to the strength of the amplitude and frequency modulation of the component. However, noise and interference between the components introduce large errors into the direct estimates, whereas the ridge estimates are more robust in the face of such effects. The best choice between the two therefore depends on the amount of noise and characteristic amplitude/frequency variation of the component to be extracted. For the time being, to avoid complications, let us assume that we reconstruct everything by a particular method, e.g.\ direct, while its adaptive choice will be discussed later in Sec.\ \ref{sec:precm}.

\subsection{Harmonics: extracting the candidates}\label{sec:hextract}

\begin{figure}[t]
\includegraphics[width=1.0\linewidth]{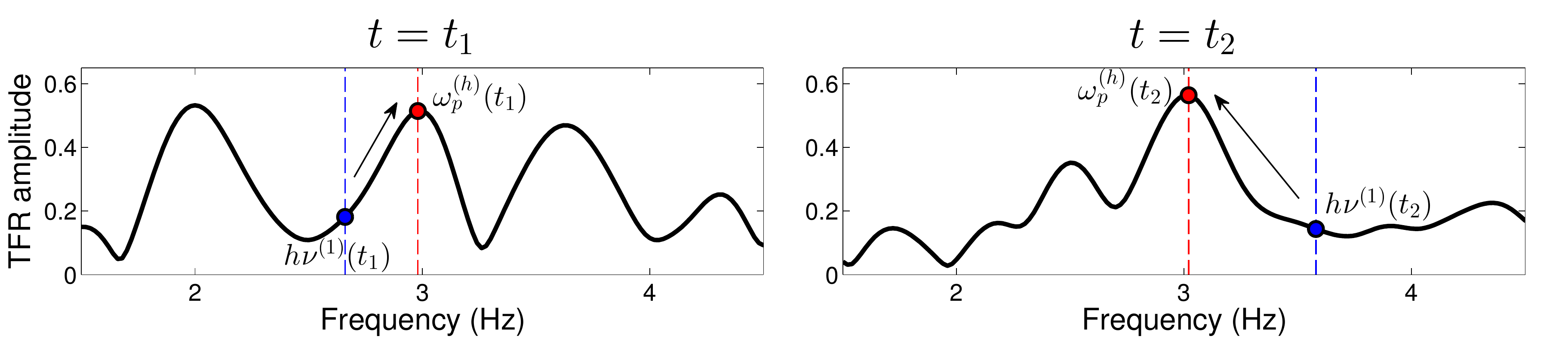}
\caption{Illustration of the extraction of the $h^{\rm th}$ harmonic ridge curve $\omega_p^{(h)}(t)$ based on the fundamental instantaneous frequency $\nu^{(1)}(t)$. At each time, starting from the expected harmonic frequency $h\nu^{(1)}(t)$ (blue points), one climbs (i.e.\ follows in the direction of TFR amplitude increase) to the nearest peak, which is then assigned to $\omega_p^{(h)}(t)$ (red points).}
\label{fig:eharm}
\end{figure}

The component extracted in the previous step will generally represent a harmonic of some Nonlinear Mode. For simplicity, however, in this and the next subsection we assume that it corresponds to a fundamental, i.e.\ the first, harmonic (we will get rid of this assumption in Sec.\ \ref{sec:hissues}). Then, given its instantaneous frequency $\nu^{(1)}(t)$, reconstructed by (\ref{ridgerec}) or (\ref{dirrec}), the harmonics for $h=2,3,...$ should be located at frequencies $h\nu^{(1)}(t)$. Therefore, the ridge curve $\omega_p^{(h)}(t)$ corresponding to the $h^{\rm th}$ harmonic can be extracted simply as the sequence of peaks which are located in the same time-frequency support (region of unimodal TFR amplitude at each time) as $h\nu^{(1)}(t)$; or, in other words, the sequence of peaks nearest to $h\nu^{(1)}(t)$ in the direction of TFR amplitude increase. This is illustrated in Fig.\ \ref{fig:eharm}. Having found $\omega_p^{(h)}(t)$, the parameters of the harmonic $A^{(h)}(t),\phi^{(h)}(t),\nu^{(h)}(t)$ can be reconstructed as discussed in the previous section.

\subsection{Harmonics: identifying the true ones}\label{sec:htest}

In general, the procedure of the previous subsection yields what is not necessarily a genuine harmonic, but just a candidate for one. Thus, even if the NM does not contain a particular harmonic, one will still obtain some signal for it, consisting of noise or components lying near the frequency of the putative harmonic. Hence, having extracted the $h^{\rm th}$ harmonic candidate, we must next determine whether it is a true harmonic or not. To tackle this problem, we use the method of surrogate data \cite{Schreiber:00,Theiler:92}, testing against the null hypothesis of independence between the first harmonic and the extracted harmonic candidate. Thus, we first select a measure to quantify the degree of dependence between the dynamics of two harmonics, and we then calculate it for the original harmonic and for a collection of surrogates -- a specially constructed time-series consistent with the null hypothesis being tested. If the original value of the selected measure lies outside the distribution of its surrogate values, this indicates genuine interdependence and the harmonic is regarded as true; otherwise, it is discarded as false.

The amplitude, phase and frequency dynamics of a true harmonic should depend on those for the first harmonic (fundamental) as $A^{(h)}(t)/A^{(1)}(t)\equiv a_h=\operatorname{const}$, $\phi^{(h)}(t)-h\phi^{(1)}(t)\equiv \varphi_h-h\varphi_1=\operatorname{const}$, $\nu^{(h)}(t)=h\nu^{(1)}(t)$. One can introduce measures $q_{A,\phi,\nu}^{(h)}\in[0,1]$ quantifying the degree of consistency with these laws ($0$ -- no consistency, $1$ - full consistency), i.e.\ the dependence of the parameters of the $h^{\rm th}$ harmonic on the corresponding parameters of the first one:
\begin{equation}\label{rhoapf}
\begin{aligned}
&q_A^{(h)}\equiv \exp\Bigg[-\frac{\sqrt{\big\langle\big[A^{(h)}(t)\langle A^{(1)}(t)\rangle-A^{(1)}(t)\langle A^{(h)}(t)\rangle\big]^2\big\rangle}}{\langle A^{(1)}(t)A^{(h)}(t)\rangle}\Bigg],\\
&q_\phi^{(h)}\equiv \big|\big\langle \exp\big[i\big(\phi^{(h)}(t)-h\phi^{(1)}(t)\big)\big] \big\rangle\big|,\\
&q_\nu^{(h)}\equiv\exp\Bigg[-\frac{\sqrt{\big\langle\big[\nu^{(h)}(t)-h\nu^{(1)}(t)\big]^2\big\rangle}}{\langle \nu^{(h)}(t)\rangle}\Bigg].\\
\end{aligned}
\end{equation}
The overall measure of interdependence between the harmonics can then be taken as
\begin{equation}\label{apc}
\rho^{(h)}(w_A,w_\phi,w_\nu)=\Big(q_A^{(h)}\Big)^{w_A}\Big(q_\phi^{(h)}\Big)^{w_\phi}\Big(q_\nu^{(h)}\Big)^{w_\nu}
\end{equation}
with the parameters $w_{A,\phi,\nu}$ giving weights to each of the consistencies $q_{A,\phi,\nu}^{(h)}$. By default, we use $\rho^{(h)}\equiv\rho^{(h)}(1,1,0)$, giving equal weights to the amplitude and phase consistencies, and no weight to the frequency consistency. The latter is because the procedure of harmonic extraction (see Sec.\ \ref{sec:hextract}), being based on the instantaneous frequency of the first harmonic, in itself introduces a dependence of $\nu^{(h)}(t)$ on $\nu^{(1)}(t)$, so it is better not to base any conclusions on it.

Ideally, for true harmonics one should have $q_{A,\phi,\nu}^{(h)}=1$, but noise, finite frequency and time resolutions, and the interference with the other components all introduce errors. In reality, therefore, the consistencies will be smaller than unity even for true harmonics. Hence one cannot identify harmonics based only on the value of $\rho^{(h)}$ but also needs to use the surrogate test. We employ the idea of time-shifted surrogates \cite{Andrzejak:03,Quiroga:02}, using as a first harmonic surrogate its time-shifted version and as the other harmonic surrogate the corresponding candidate harmonic extracted from the time-shifted TFR. Such time-delay destroys any temporal correlations between the signals while preserving all other their properties, thus creating surrogates consistent with the null hypothesis of independence.

Given the maximal time-shift (in samples) $M$, the surrogate parameters for the first harmonic are taken as shifted on $\Delta T_d/2$ backwards
\begin{equation}\label{tsurr1}
\begin{aligned}
&A_d^{(1)}(\tau)=A^{(1)}(\tau-\Delta T_d/2),\;\phi_d^{(1)}(\tau)=\phi^{(1)}(\tau-\Delta T_d/2),\;\nu_d^{(1)}(\tau)=\nu^{(1)}(\tau-\Delta T_d/2),\\
&\tau=\{t_{i=1+M/2,\dots,N-M/2}\},\;\Delta T_{d=1,\dots,N_d}=M(1-2d/N_d)/2f_s,
\end{aligned}
\end{equation}
where $f_s$ is the signal sampling frequency and $N$ is its total length in samples (note, that the length of the surrogate time series is smaller than the original ones, being $N-M$). Using $\nu_d^{(1)}(\tau)$ (\ref{tsurr1}) as a reference profile, the surrogate $h^{\rm th}$ harmonic and its parameters $A^{(h)}_d(\tau),\phi^{(h)}_d(\tau),\nu^{(h)}_d(\tau)$ are extracted in the same way as described in Sec.\ \ref{sec:hextract}, but from the signal's TFR shifted on $\Delta T_d/2$ forward in time ($G(\omega,\tau+\Delta T_d/2)$ or $W(\omega,\tau+\Delta T_d/2)$).

The extraction of the corresponding curves for all surrogates can greatly be accelerated by initial preprocessing of the TFR, constructing its skeleton \cite{Iatsenko:13tfr2} at the beginning and utilizing it for each surrogate. Thus, at each time $t$ we break the TFR into the regions of unimodal TFR amplitude $[\omega_{-,m}(t),\omega_{+,m}(t)]$, $m=1,2,\dots,N_p(t)$, and from them reconstruct by (\ref{ridgerec}) or (\ref{dirrec}) the corresponding amplitudes $\overline{A}_m(t)$, phases $\overline{\phi}_m(t)$ and frequencies $\overline{\nu}_m(t)$ (which form the TFR skeleton). Then what is left for each surrogate is to find the consequence of the indices $m_d(t)$ of the supports in which the expected harmonic frequency lies and to take the corresponding parameters:
\begin{equation}\label{tsurr2}
\begin{aligned}
&m_d(\tau):\;h\nu_d^{(1)}(\tau)\equiv h\nu^{(1)}(\tau-\Delta T_d/2)\in[\omega_{-,m_d(\tau)}(\tau+\Delta T_d/2),\omega_{+,m_d(\tau)}(\tau+\Delta T_d/2)]\\
&A_d^{(h)}(\tau)=\overline{A}_{m_d(\tau)}(\tau),\;\phi_d^{(h)}(\tau)=\overline{\phi}_{m_d(\tau)}(\tau),\;\nu_d^{(h)}(\tau)=\overline{\nu}_{m_d(\tau)}(\tau)\\
\end{aligned}
\end{equation}

\begin{remark}
Alternatively, one can construct the surrogates by shifting the already reconstructed parameters of the first and $h^{\rm th}$ harmonic relative to each other, without re-extracting them for each surrogate. In fact, we have found that such an approach gives almost the same results and at the same time is faster. However, as mentioned before, the original procedure of harmonic extraction introduces the correlation between $\nu^{(h)}(t)$ and $\nu^{(1)}(t)$, and can therefore introduce dependence between other parameters, such as phases. Therefore, to be rigorous, the procedure should be repeated for each surrogate, thus automatically taking into account possible bias of this kind and producing more reliable estimates of the significance.
\end{remark}

Summarizing, we calculate the amplitude-phase consistencies $\rho_{d=1,\dots,N_d}^{(h)}(1,1,0)$ (\ref{apc}) for the surrogate parameters (\ref{tsurr1}), (\ref{tsurr2}) and compare them with the value $\rho_{0}^{(h)}(1,1,0)$ calculated in the same way but for the zero time shift $\Delta T_0=0$. The probability measure (although mathematically not the true probability) that the extracted $h^{\rm th}$ harmonic curve is a true harmonic of the main one is then quantified by the significance of the surrogate test, i.e.\ by the relative part of surrogates for which $\rho_d^{(h)}<\rho_0^{(h)}$. For example, if we found $1000$ surrogate amplitude-phase consistencies $\rho_{d=1,..,1000}^{(h)}$ and $792$ of them are smaller than the original value $\rho_0^{(h)}$, then the rough probability that the extracted curve represents a true harmonic is $79.2\%$. We regard a harmonic as true if the probability calculated in this way is $\geq 95\%$. By default we use $N_d=100$ surrogates and a maximum time shift that equals quarter-length of the signal $M=N/4$, so that the surrogates are of $N-M=3N/4$ length each.

Note that the significance of the surrogate test does not depend on the magnitude of $\rho^{(h)}$ (\ref{apc}). Thus, there might be an independent component located at the frequency of a possible harmonic, so that the amplitude-phase consistency would be high but, because it does not fully adjust its amplitude and phase to that of the fundamental harmonic, the surrogate test will not reject the null hypothesis of independence (e.g.\ see below, Table 1). Thus, the possibility of picking up a spurious component that is nearby in frequency is largely eliminated.

\begin{remark}
Since we test many harmonic candidates, it is natural to expect that sometimes we will encounter false positives, i.e.\ the surrogate test will regard as true a harmonic which is actually false. Thus, for some noise realizations the parameters of the TFR around the expected harmonic frequency might indeed appear to be correlated with the parameters of the fundamental harmonic. However, in such cases the extracted harmonic candidate will usually have quite small amplitude-phase consistency $\rho^{(h)}$ (\ref{apc}). To reduce the probability of false positives, therefore, we introduce some threshold $\rho_{\min}$ and regard the harmonic as true only if it both passes the surrogate test and at the same time is characterized by $\rho^{(h)}\geq\rho_{\min}$. Empirically, we set this threshold as
\begin{equation}\label{rhomin}
\rho_{\min}=0.5^{w_A+w_\phi},
\end{equation}
where $w_A,w_\phi$ are the weightings used in (\ref{apc}) (as mentioned, we use $w_A=w_\phi=1$); the value (\ref{rhomin}) was found to be quite universal and to work well in most cases. Note, that for a true harmonic one can also have $\rho^{(h)}<\rho_{\min}$, but this will usually mean that it is badly damaged by noise or other influences and cannot be recovered without large errors.
\end{remark}

\subsection{Harmonics: practical issues}\label{sec:hissues}

\textbf{Extracting in order.} To improve the accuracy of reconstruction, each harmonic which is identified as true should be subtracted from the signal before extracting and testing the next one. This decreases the errors related to interference between the harmonics and makes all procedures more accurate. The same consideration applies to the first harmonic which, after being found by the methods described in Sec.\ \ref{sec:cextract}, should be subtracted from the signal before extraction of any of the other harmonics.

\textbf{How many harmonics to extract?} Clearly, the maximum number of harmonics one can extract in principle is $h_{\max}=(f_s/2)/\langle\nu^{(1)}(t)/2\pi\rangle$, where $f_s$ is the sampling frequency of the signal (so that $f_s/2$ is the maximally achievable, Nyquist frequency) and $\nu^{(1)}(t)$ denote the extracted instantaneous frequency of the first harmonic. However, checking all harmonics might be computationally very expensive, and is often not needed. In practice, it is better to stop the search after some chosen number of sequential harmonics has been identified as false, making it likely that there are no more true harmonics. We choose this number to be $S=3$ by default.

\textbf{What if the extracted component is not the first harmonic?} Although intuitively the first harmonic should have the highest amplitude (and will therefore be extracted as the dominant curve), for some complicated waveforms (e.g.\ with sharp peaks) this might be untrue. Therefore, before starting extraction of the harmonics, one should first ensure that one is starting with the first harmonic; and, if not, then find it and switch to it. To do this, we apply the same procedure described for harmonic extraction and identification, but in the reverse direction, i.e.\ using $h=1/2,1/3,1/4,...$. Then if some of these are identified as true, we switch to the one with the smallest frequency and start from there. The minimal $h$ one can go for can be set as $h_{\min}=(1/T)/\langle\nu^{(1)}(t)/2\pi\rangle$, although the statistics will be bad already for $h<5h_{\min}$, when the related oscillation has less than 5 cycles over the whole signal time-duration \cite{Iatsenko:13tfr1}. Nevertheless, it is better to stop after $S_0=3$ consecutive $1/n$ harmonics have been identified as false, in the same manner as is done for the usual harmonics.

\subsection{Stopping criterion}\label{sec:stopcrit}

Once all harmonics are identified and reconstructed, they are summed up into the Nonlinear Mode, which is then subtracted from the signal and the procedure is repeated on the residual. The natural question arises, therefore, of how many Nonlinear Modes to extract, i.e.\ when to stop the decomposition. Obviously, decomposition of any noise, e.g.\ white, Brownian, any other correlated or not, does not make much sense, so the reasonable goal is to extract all oscillatory components present in the signal and leave the noise and trends as the residual. Therefore, after extraction of each NM one needs to decide whether what is left contains any more meaningful oscillations, in which case one continues the decomposition, or whether it just represents noise, in which case one should stop. The problem is thus reduced to distinguishing deterministic from random dynamics.

To solve it, one can use the surrogate test against the null hypothesis of linear noise \cite{Schreiber:00,Theiler:92}, which includes white and colored noises (e.g. Brownian). The surrogates for this task, called FT surrogates, are constructed by taking the inverse Fourier transform of the signal's FT with randomized phases of the Fourier coefficients: $s_s(t)=(2\pi)^{-1}\int [\hat{s}(\xi)e^{i\varphi_s(\xi)}]e^{i\xi t}d\xi$, where $\varphi_s(\xi)=-\varphi_s(-\xi)$ denote the phases taken at random uniformly on $[0,2\pi]$ for each frequency $\xi>0$. The reason for this is that any linear noise (an ARMA process $x(t_n)=a_0+b_0\eta_W(t_n)+\sum_{p=1}^Pa_px(t_{n-p})+\sum_{m=1}^Mb_m\eta_W(t_{n-p})$, where $\eta_W(t_n)$ denotes Gaussian white noise of unit variance) is characterized only by the amplitudes of the Fourier coefficients. Randomization of the Fourier phases preserves the power spectrum, so that the surrogate time series will represent another realization of the same random process if the original time series is noise, thus being consistent with the tested null hypothesis. On the other hand, if some meaningful NMs are present in the signal, the randomization of the phases will destroy the particular phase relationships responsible for the amplitude and frequency modulations, making their behavior less deterministic.

\begin{remark}
Instead of this stopping criterion, one can use the conventional approach of stopping when the standard deviation of the residual is lower than some percentage of that for the original signal. Another way is to extract some predefined number of NMs, e.g.\ 10. However, the criterion based on the surrogate test seems to be the most physically meaningful.
\end{remark}

One now needs to select the discriminating statistics, which is calculated for the original signal and the surrogates, so that the null hypothesis of linear noise is accepted if the original value lies within the surrogate values and rejected otherwise. The commonly used statistics involve one of the generalized Renyi dimensions \cite{Hentschel:83,Grassberger:83,Farmer:83}, with the correlation dimension calculated by Grassberger-Procaccia approach \cite{Grassberger:83cd1,Grassberger:83cd2,Theiler:87} remaining the most popular choice. However, we find that the surrogate test based on such measures is extremely sensitive to noise, prohibiting their use in the NMD which is intended to be noise-robust.

Therefore, we need to devise another discriminating statistics. There are virtually no restrictions on its choice \cite{Theiler:96}, so that almost any measure can be used. The only question to be asked is how powerful it is, i.e.\ how good in distinguishing the deterministic dynamics from noise. Given that NMD is based on the WFT/WT, it seems reasonable to select statistics based on the properties of the time-frequency representation obtained. Namely, since at the first step of the NMD process we extract the component from the signal's TFR (see Sec.\ \ref{sec:cextract}), we can base our statistics on the properties of the components extracted (in the same way) from the original signal and from its FT surrogates. Thus, if the original component is true (and not just formed from noise peaks picked in the time-frequency plane), then it is expected to have more deterministic amplitude and frequency modulation than the surrogate components, which should be more stochastic; otherwise, there will be no difference.

The degree of order in the extracted amplitude $A(t)$ or frequency $\nu(t)$ can be quantified by their spectral entropies $Q[\hat{A}(\xi)],Q[\hat{\nu}(\xi)]$, so that the discriminating statistics $D$ for the surrogate test can be taken as their combination
\begin{equation}\label{sentropy}
\begin{aligned}
&D(\alpha_A,\alpha_\nu)\equiv \alpha_A Q[\hat{A}(\xi)]+\alpha_\nu Q[\hat{\nu}(\xi)],\\
&Q[f(x)]\equiv-\int \frac{|f(x)|^2}{\int |f(x)|^2dx}\log \frac{|f(x)|^2}{\int |f(x)|^2dx}dx.
\end{aligned}
\end{equation}
Note that, in practice, due to the finite sampling frequency $f_s$ and sample length $N$ of real signals, the integrals over frequency $\xi$ in $Q[\hat{A}(\xi)],Q[\hat{\nu}(\xi)]$ (\ref{sentropy}) are discretized into sums over the discrete FT frequencies $\xi_n=(n/N-1/2)f_s$, $n=1,\dots,N$.

In the present context, we have found the statistics $D(\alpha_A,\alpha_\nu)$ (\ref{sentropy}) to be more meaningful and much more powerful than other choices (e.g.\ the popular correlation dimension \cite{Grassberger:83cd1,Grassberger:83cd2,Theiler:87}). This statistics is directly related to the quality of the representation of component in the signal's TFR, so that the significance of the surrogate test based on it reflects the proportion of the ``deterministic'' part in the extracted amplitude and frequency dynamics. Thus, if the residual does not pass the surrogate test (null hypothesis is not rejected), this might mean either that the residual is indeed noise, or that the component simply cannot be reliably extracted from the TFR (e.g.\ because resolution characteristics of the TFR are not appropriate to reliably represent the related amplitude/frequency modulation and/or to segregate the component from noise).

The power of $D(\alpha_A,\alpha_\nu)$, i.e.\ its ability to distinguish deterministic from random dynamics, depends strongly on the complexity of the component's amplitude and frequency modulations: the lower the original spectral entropies $Q[\hat{A}(\xi)],Q[\hat{\nu}(\xi)]$ are, the more powerful is the test. However, even in the (quite unrealistic) case when the signal contains a meaningful component without any amplitude or frequency modulation, i.e.\ a pure tone $A\cos \nu t$, due to numerical issues \cite{Theiler:93,Stam:98} the surrogate test will still be quite powerful in rejecting the null hypothesis (unless this tone has an integer number of cycles over the time-length of the signal). The power of the test is also inversely proportional to the spread of $\hat{A}(\xi),\hat{\nu}(\xi)$: the more concentrated they are, the narrower the frequency band that the component $A(t)\cos \phi(t)$ occupies, so that the less the noise power that is contained in it.

As to the choice of $\alpha_A,\alpha_\nu$ in (\ref{sentropy}), we have found that the powers of $D(1,0)=Q[\hat{A}(\xi)]$ and $D(0,1)=Q[\hat{\nu}(\xi)]$ are often inversely proportional to the strengths of the amplitude and frequency modulation, respectively. Thus, $D(1,0)$ is preferable for components with relatively small amplitude variability and considerable frequency variability, while $D(0,1)$ is better otherwise. Although the procedure is not fully rigorous  mathematically, we therefore perform three tests, using $D(1,0)$, $D(0,1)$ and $D(1,1)$, and then select the significance as the maximum among them. It remains to be established, however, whether some better statistics not having the drawback mentioned can be found.

Summarizing, we extract the components from the TFR of the original signal and compute the corresponding $D_0(\alpha_A,\alpha_\nu)$ (\ref{sentropy}); then we create $N_s$ FT surrogates of the signal, for each of them calculate the corresponding TFR, extract the component from it and compute the respective $D_{s=1,\dots,N_s}(\alpha_A,\alpha_\nu)$. We use $N_s=40$ surrogates and set the significance level to $95\%$, rejecting the tested null hypothesis of noise if the number of surrogates with $D_s>D_0$ is equal or higher than $0.95\times40=38$. The test is performed for three different $(\alpha_A,\alpha_\nu)$ in (\ref{sentropy}), using $D(1,1)$, $D(0,1)$ and $D(1,0)$ as a discriminating statistic; if at least for one of them the null hypothesis is rejected, we regard the signal as not noise and continue the decomposition.


\begin{remark}
At the very start of each NMD iteration, the signal should be detrended, which we do by subtracting a third order polynomial fit from it. This is especially important for the surrogate test, as trends might greatly affect its performance, In particular, the FT surrogates should be constructed using an already detrended signal as a base, because otherwise one might end up with testing the trend against noise rather than the oscillatory component of interest. Additionally, to guarantee that the boundary distortions are of the same nature in the original TFR and those for the surrogates, it is recommended to use padding with zeros \cite{Iatsenko:13tfr1}; the original $D_0$ should thus be recalculated using parameters of the component extracted from the TFR of the zero-padded signal. This issue, however, concerns only the surrogate test, while in all other cases we use the more accurate (but at the same time more computationally expensive) predictive padding \cite{Iatsenko:13tfr1}.
\end{remark}

\begin{remark}
Instead of the FT surrogates used here, one can alternatively utilize the more advanced IAAFT surrogates \cite{Schreiber:96}, which are devoted to testing against the null hypothesis of an invertibly rescaled linear stochastic process (e.g.\ Brownian noise taken to the third power). These surrogates preserve both the power spectrum and, as much as possible, the distribution of values of the original signal. However, as can be seen from (\ref{wftwt}), the WFT and WT do not explicitly take into account the distribution of signal values, only its FT. Therefore, the simpler and easier-to-compute FT surrogates seem to be a more natural choice in our case, and we have found that the test developed above works well even for rescaled (either invertibly or non-invertibly), non-Gaussian and non-stationary linear noise. In fact, the outcome of the test seems to be determined primarily by the ``goodness'' of the component representation in TFR, which is affected only by the noise power and its distribution in the time-frequency region around the component's instantaneous frequency.
\end{remark}

\section{Improvements}\label{sec:improve}

The NMD as described in the previous section already represents a very powerful decomposition tool. However, it may be made even better with the help of the upgrades outlined below, although at the expense of greatly increased computational cost for some of them (though it still remains $O(N\log N)$).

\subsection{Adaptive representation of the harmonics}\label{sec:ihadapt}

Even if the first harmonic of some NM is well represented in the TFR and can be accurately extracted and reconstructed, it does not mean that the same applies to all the other harmonics too. For example, harmonics of an NM with only amplitude modulation require the same time and frequency resolution, so that the WFT can represent them all well, while for the WT, where time and frequency resolution scale with frequency, one will need to adjust the resolution parameter $f_0$ for each harmonic. Thus, consider an NM with simple sinusoidal amplitude modulation:
\begin{equation}\label{NMam}
\begin{aligned}
s(t)=&(1+r_a\cos(\nu_a t+\varphi_a))\sum_{h=1}^\infty a_h\cos (h\nu t+\varphi_h)\equiv\sum_{h=1}^\infty x^{(h)}(t)\\
\Rightarrow x^{(h)}(t)&\equiv a_h(1+r_a\cos(\nu_a t+\varphi_a))\cos (h\phi(t)+\varphi_h),\;\phi(t)=\nu t\\
&=a_h\big[\cos (h\nu t+\varphi_h)+\frac{r_a}{2}\cos((h\nu+\nu_a) t+(\varphi_h+\varphi_a))+\frac{r_a}{2}\cos((h\nu-\nu_a) t+(\varphi_h-\varphi_a))\big].
\end{aligned}
\end{equation}
Clearly, all harmonics have the same amplitude variability (in relative terms), so that the TFR time-resolution should also be the same for them. This is reflected in the fact that each harmonic $x^{(h)}(t)$ is composed of three Fourier terms, which all have identical amplitude ratios, frequency differences and phase relationships for each $h$. Furthermore, the frequency distance between two consecutive harmonics remains the same, meaning that the frequency resolution also should not be changed for harmonics.

Therefore, the WFT, having constant time and frequency resolution, will be a perfect match for this case. This means, first, that if the amplitude modulation of the first harmonic is represented reliably in the WFT, then the same will also apply to all other harmonics and, secondly, that if two first harmonics do not interfere in the WFT, then any two harmonics will also be well-separated. For the WT, the former is also true, as the time-resolution, i.e.\ the ability to reflect changes in time, increases with frequency for this type of TFR. However, the frequency resolution of the WT progressively worsens with the increase of frequency, so that the higher are the harmonics, the harder it is to resolve them.

\begin{figure}[t]
\includegraphics[width=1.0\linewidth]{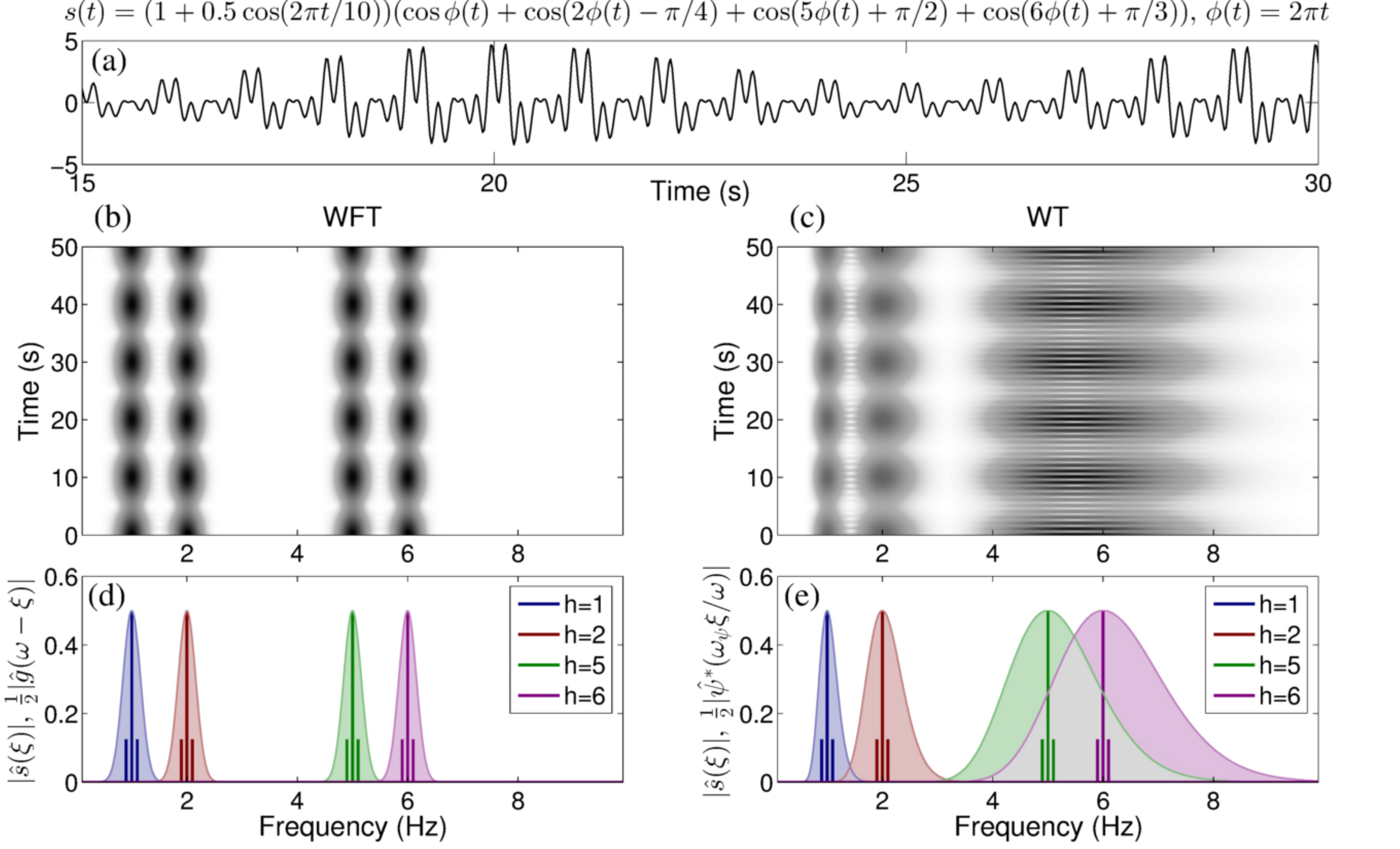}
\caption{(a) Nonlinear Mode with amplitude modulation, as specified by the equation for $s(t)$ above the box. (b),(c) The corresponding WFT and WT amplitudes, respectively. (d),(e) The signal's FT, integrated over one frequency bin (vertical lines), with the parts responsible for different harmonics shown in different colors; the shaded areas show the absolute values of the window functions $\hat{g}(\omega-\xi)$ or wavelet functions $\hat{\psi}^*(\omega_\psi\xi/\omega)$ centered at the mean frequencies of the harmonics $\omega=2\pi h$ and rescaled to half of the harmonics' mean amplitudes. The signal is sampled at 50 Hz for 50 s.}
\label{fig:harm1}
\end{figure}

This is illustrated in Fig.\ \ref{fig:harm1}, where from (b) and (c) it is clear that all harmonics can be well represented in the WFT, but for the WT higher harmonics begin to interfere. This issue is explained schematically in Fig.\ \ref{fig:harm1}(d) and (e). Thus, both the WFT and WT can be seen as convolutions in the frequency domain of the signal with a window $\hat{g}(\omega-\xi)$ and wavelet $\hat{\psi}^*(\omega_\psi\xi/\omega)$, as seen from (\ref{wftwt}). Figs.\ \ref{fig:harm1}(d) and (e) show the signal's discrete FT $\hat{s}(\xi)$ together with the (rescaled and centered at the mean frequencies of the harmonics $\omega/2\pi=1,2,5,6$ Hz) window and wavelet FTs $\hat{g}(\omega-\xi)$ and $\hat{\psi}^*(\omega_\psi\xi/\omega)$ (\ref{winwav}), with which $\hat{s}(\xi)$ is convoluted while constructing the WFT and WT, respectively.

Roughly speaking, the quality of the representation of time-variability (the amplitude modulation in the present case) for each harmonic can be estimated based from the proportion of the corresponding Fourier coefficients lying below the shown window/wavelet FTs. For both the WFT and WT, each harmonic has three Fourier coefficients (\ref{NMam}), all of which lie in the correspondingly shaded areas, meaning that the time-resolution in each case is sufficient to represent amplitude modulation appropriately. The degree of interference between the harmonics in the WFT and WT can be estimated from the area of overlap between the corresponding window/wavelet FTs. For the WFT (Fig.\ \ref{fig:harm1}(d)) the interference between harmonics does not depend on their frequency, while for the WT (Fig.\ \ref{fig:harm1}(e)) it increases for the higher harmonics, which results in the fifth and sixth harmonics being represented as a single component (Fig.\ \ref{fig:harm1}(c)). This is because the minimal frequency difference between two harmonics is equal to the frequency of the first harmonic, thus being defined on a linear scale (which is natural for the WFT), but not on a logarithmic one (natural for the WT).

The situation changes when there is frequency modulation. In this case neither the WFT nor the WT provide optimal representation of the harmonics; in addition, some harmonics cannot in principle be reliably resolved with time-frequency analysis methods. Thus, consider an NM with simple sinusoidal frequency modulation
\begin{equation}\label{NMfm}
\begin{aligned}
s(t)=&\sum_{h=1}^\infty a_h\cos (h\phi(t)+\varphi_h)\equiv\sum_{h=1}^\infty x^{(h)}(t),\;\phi(t)=\nu t+r_b\sin(\nu_b t+\varphi_b)\\
\Rightarrow x^{(h)}(t)&=a_h\cos (h\nu t+hr_b\sin(\nu_b t+\varphi_b)+\varphi_h)
=a_h\operatorname{Re}\Big[\sum_{n=-\infty}^\infty J_n(hr_b)e^{i(\varphi_h+n\varphi_b)}e^{i(h\nu+n\nu_b)t}\Big]
\end{aligned}
\end{equation}
where we have used the formula $e^{ia\sin\phi}=\sum_{n=-\infty}^{\infty}J_n(a)e^{in\phi}$ with $J_n(a)=(-1)^nJ_{-n}(a)$ denoting a Bessel function of the first kind. Since after some $|n|$ all terms in (\ref{NMfm}) become negligible, one can in practice restrict the summation to $|n|\leq n_J(ha)$, with the maximum non-negligible order $n_J$ being determined as
\begin{equation}\label{nj}
n_J(a):\;\frac{\sum_{|n|>n_J(a)}|J_n(a)|^2}{\sum_{n=-\infty}^\infty|J_n(a)|^2}<\epsilon_p,
\end{equation}
where $\epsilon_p$ denote a chosen accuracy threshold. The value of $n_J(a)$ increases with $a$, being (for $\epsilon_p=0.001$): $n_J(0\lesssim a\lesssim 0.3)=2$, $n_J(0.3\lesssim a\lesssim 0.65)=3$, $n_J(0.65\lesssim a\lesssim 1.13)=4$ etc. As a result, the higher the harmonic, the larger the frequency range it occupies, i.e.\ the larger is the number of non-negligible terms in (\ref{NMfm}).

\begin{figure}[t]
\includegraphics[width=1.0\linewidth]{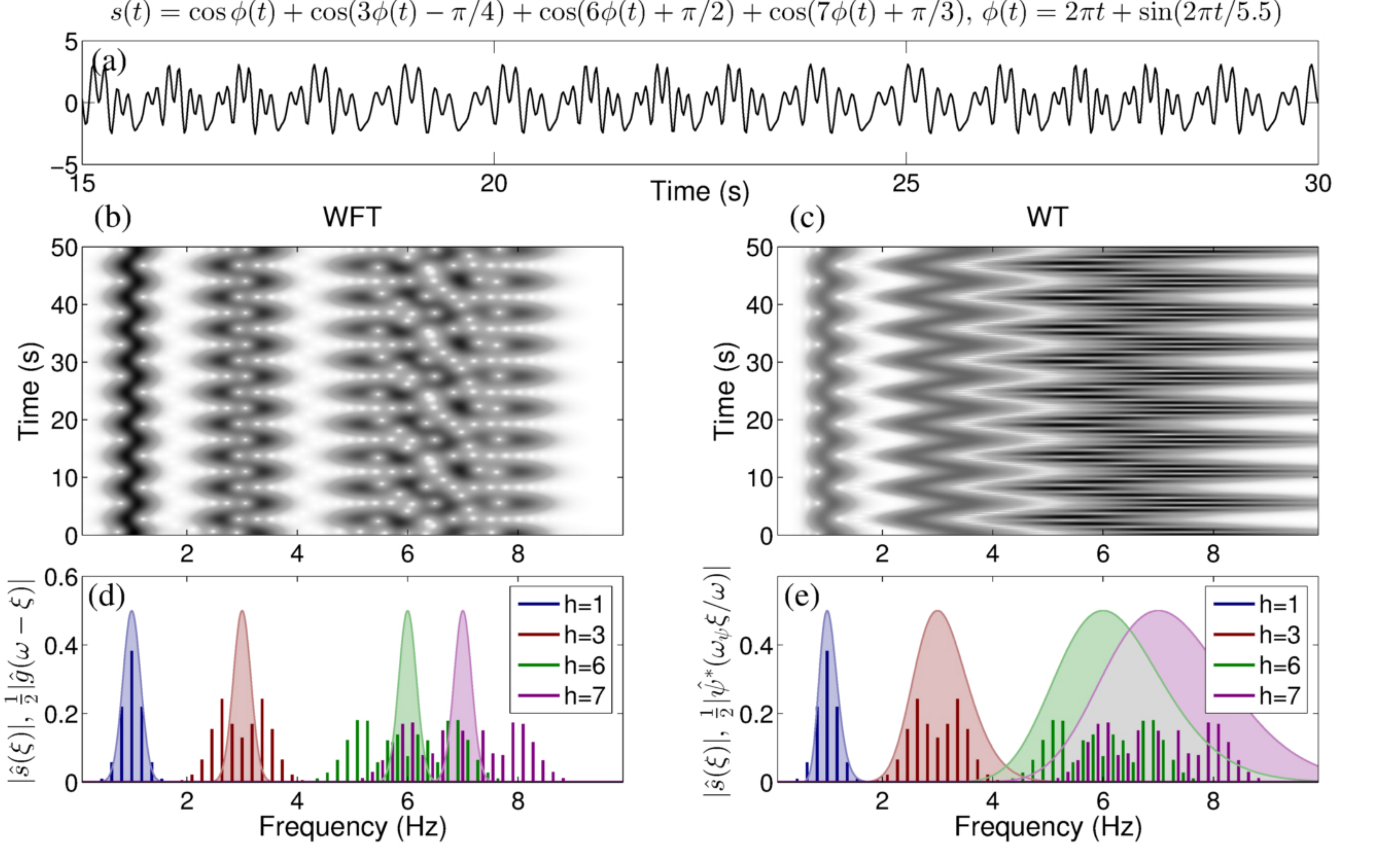}
\caption{Same as in Fig.\ \ref{fig:harm1}, but for the frequency-modulated NM plotted in (a) and specified by the equation for $s(t)$ above the box. The signal is sampled at 50 Hz for 55 s.}
\label{fig:harm2}
\end{figure}

Consequently, to reflect reliably the frequency modulation of the higher harmonics, one needs higher time-resolution for them, a requirement that is satisfied by the WT, but not by the WFT. This issue is illustrated in Fig.\ \ref{fig:harm2}, where it can be seen that the WT can represent reliably both the first and third harmonics, while the WFT reflects appropriately only the first one. However, the increased time-resolution of the WT is provided at the expense of decreased frequency resolution, leading to stronger interference between harmonics, as seen from the previous case (Fig.\ \ref{fig:harm1}). Figure \ref{fig:harm2} also shows, that in some cases it might in principle be impossible to represent reliably two harmonics in the TFR. Thus, as can be seen from Fig.\ \ref{fig:harm2}(d) and (e), the FTs $\hat{x}^{(h)}(\xi)$ of the sixth and seventh harmonics are ``entangled'', i.e.\ the frequency regions in which they are contained overlap. Therefore, unless specifically designed for this case, any window/wavelet function which picks the Fourier coefficients of one harmonic will inevitably also pick those corresponding to the other one too.

\begin{remark}\label{rem:emdres}
Interestingly, the (E)EMD \cite{Huang:98,Wu:09} has logarithmic frequency resolution \cite{Feldman:09,Rilling:08}, and therefore suffers from the same drawbacks as the WT in relation to representation of harmonics. Thus, (E)EMD applied to the NM with amplitude modulation shown in Fig.\ \ref{fig:harm1} will decompose it into its first and second harmonics, but will merge the fifth and sixth ones into a single mode. Note also that, while for the WT the logarithmic frequency resolution can be tuned by changing the resolution parameter $f_0$ in (\ref{wftwt}), for (E)EMD the resolution is fixed around some implicit value \cite{Feldman:09,Rilling:08}.
\end{remark}

Summarizing, for the general case where both amplitude and frequency modulation are present in the NM, accurate representation of higher harmonics requires higher time resolution, but the same frequency resolution. However, an increase in time resolution inevitably leads to a decrease of frequency resolution, so usually one needs to search for some compromise. Due to this, it is often the case that neither the WFT nor the WT with a constant resolution parameter $f_0$ in (\ref{winwav}) can represent all harmonics reliably.

To tackle this problem, one can adaptively adjust $f_0$ for each harmonic individually. Assuming that the extracted first harmonic is reconstructed well, the quality of the representation of the $h^{\rm th}$ harmonic can be assessed through its consistency $\rho^{(h)}$ (\ref{apc}) with the first one. Therefore, the optimal resolution parameter for the $h^{\rm th}$ harmonic, $f_0^{(h)}$, can be selected as that for which $\rho^{(h)}$ (\ref{apc}) achieves its maximum and at the same time the harmonic passes the surrogate test, i.e.\ is identified as a true harmonic.

The remaining question is the region within which to search for $f_0^{(h)}$. To find it, we first need to understand in what frequency band $[\omega_f^{(h)}-\Delta\omega_f^{(h)}/2,\omega_f^{(h)}+\Delta\omega_f^{(h)}/2]$ the $h^{\rm th}$ harmonic FT $\hat{x}^{(h)}(\xi)$ is concentrated, given the frequency band of the first harmonic. We define $\omega_f^{(h)}$ and $\Delta\omega_f^{(h)}$ as
\begin{equation}\label{fband}
\omega_f^{(h)},\Delta\omega_f^{(h)}:\;
\int_0^{\omega_f^{(h)}-\Delta\omega_f^{(h)}/2}|\hat{x}^{(h)}(\xi)|^2d\xi=\frac{\epsilon_p}{2} \frac{E^{(h)}_{tot}}{2},\;\;
\int_{\omega_f^{(h)}+\Delta\omega_f^{(h)}/2}^\infty|\hat{x}^{(h)}(\xi)|^2d\xi=\frac{\epsilon_p}{2} \frac{E^{(h)}_{tot}}{2},
\end{equation}
where $E^{(h)}_{tot}\equiv \int |\hat{x}^{(h)}(\xi)|^2d\xi=\int |x^{(h)}(t)|^2dt$ is the total energy of the harmonic and $\epsilon_p$ denotes the chosen accuracy threshold (note that, as everywhere else, we assume that the analytic approximation (\ref{aa}) holds for $x^{(h)}(t)=A^{(h)}(t)\cos \phi^{(h)}(t)$). Then it can be shown that, given $\omega_f^{(1)},\Delta\omega_f^{(1)}$, the frequency range for the $h^{\rm th}$ harmonic will be
\begin{equation}\label{hband}
[\omega_f^{(h)}-\Delta\omega_f^{(h)}/2,\omega_f^{(h)}+\Delta\omega_f^{(h)}/2]:\;
\omega_f^{(h)}=h\omega_f^{(1)},\;\Delta\omega_f^{(1)}\lesssim\Delta\omega_f^{(h)}\lesssim h\Delta\omega_f^{(1)}
\end{equation}

It is clear from the preceding discussion that, for the simple AM mode (\ref{NMam}), one has $\omega_f^{(h)}=h\omega_f^{(1)}$, $\Delta\omega_f^{(h)}=\Delta\omega_f^{(1)}$. The same is true for any AM mode, as can be shown by representing the amplitude modulation through the Fourier series $A(t)=\sum_nr_n\cos(\nu_nt+\varphi_n)$ and expanding each harmonic in the same way as in (\ref{NMam}). For the FM mode (\ref{NMfm}), the expression (\ref{hband}) follows from the fact that, at least for small enough $\epsilon_p$ (e.g.\ 0.001) in the definition (\ref{nj}), one has
\begin{equation}\label{njh}
n_J(a)\lesssim n_J(h a)\lesssim hn_J(a),
\end{equation}
for most $a$, as can be confirmed numerically. Then, taking into account (\ref{njh}) in the expansion (\ref{NMfm}), one obtains (\ref{hband}). For the $h^{\rm th}$ harmonic of the NM with a more complicated double-sinusoidal frequency modulation one has
\begin{equation}\label{hfmts}
\begin{aligned}
x^{(h)}(t)=&A\cos(h\nu t+\varphi+hr_1\sin(\nu_1 t+\varphi_1)+hr_2\sin(\nu_2 t+\varphi_2))\\
=&A\operatorname{Re}\Big[e^{i(h\nu t+\varphi)}\sum_{n=-\infty}^\infty J_n(hr_1)e^{i(n\nu_1 t+n\varphi_1)}
\sum_{m=-\infty}^\infty J_m(hr_2)e^{i(m\nu_2 t+m\varphi_2)}\Big]\\
=&\sum_{n=-\infty}^\infty AJ_n(hr_1)\operatorname{Re}\Big[e^{i(h\nu t+n\nu_1 t+\varphi+n\varphi_1)}
\sum_{m=-\infty}^\infty J_m(hr_2)e^{i(m\nu_2 t+m\varphi_2)}\Big],
\end{aligned}
\end{equation}
Clearly, it can be viewed as a sum of components with singe-sinusoidal frequency modulation, but amplitudes $AJ_n(hr_1)$ and frequencies $(h\nu+n\nu_1)$. The frequency range of each of these components, and the number of them with a non-negligible amplitudes, both scale with harmonic number $h$ according to (\ref{hband}), (\ref{njh}). Therefore, the frequency range of the whole $x^{(h)}(t)$ (\ref{hfmts}) satisfies (\ref{hband}). Since frequency modulation can always be expanded in a Fourier series, using the same trick as with (\ref{hfmts}), it can be shown by induction that (\ref{hband}) holds for any FM mode. The generalization to the case of any joint amplitude and frequency modulation is then straightforward.

Based on (\ref{hband}) and the scaling properties of the WFT/WT, one can now determine the region within which the optimal resolution parameter $f_0^{(h)}$ for the $h^{\rm th}$ harmonic lies. Given that the first harmonic was extracted from the TFR calculated with resolution parameter $f_0^{(1)}$, and assuming that the corresponding time and frequency resolutions are appropriate for this first harmonic, one has
\begin{equation}\label{optfh}
\mbox{\textbf{WFT:}}\;f_0^{(h)}\in[f_0^{(1)}/h,f_0^{(1)}],\;\;\;\mbox{\textbf{WT:}}\;f_0^{(h)}\in[f_0^{(1)},hf_0^{(1)}]
\end{equation}

We search for the optimal $f_0^{(h)}$ by first breaking the region (\ref{optfh}) into $N_r$ values $f_{r=1,\dots,N_r}^{(h)}$ (we use $N_r=10$). For each of them we calculate the TFR, extract the harmonic from it (see Sec.\ \ref{sec:hextract}), estimate the corresponding consistency $\rho_r^{(h)}$ (\ref{apc}), and test the harmonic for being true (see Sec.\ \ref{sec:htest}). Among the $f_r^{(h)}$ for which the harmonic was identified as true, we select the one characterized by the highest $\rho_r^{(h)}$. It is then used as the starting guess for an iterative golden section search of the optimal $f_0^{(h)}$ (with default accuracy being $\epsilon_f=0.01$), maximizing the consistency $\rho^{(h)}$ (\ref{apc}).

Note that (\ref{optfh}) does not take into account the interference between harmonics and with the other components which might lie nearby in frequency, so that in some cases the upper bound on $f_0^{(h)}$ might be higher than (\ref{optfh}); the same consideration applies to the lower bound. Therefore, if near the minimum or maximum of the tested $f_0^{(h)}$ the consistency $\rho^{(h)}$ is found to grow, we continue the search in that direction until the peak appears.

\begin{remark}
We perform an identical procedure when we check whether the extracted component represents a first or a higher harmonic by extracting and testing its $h=1/2,1/3,...$ harmonic candidates (see Sec.\ \ref{sec:hissues}). The formulas (\ref{hband}) and (\ref{optfh}) remain valid in this case, but the upper and lower bounds for $\Delta\omega_f^{(h)}$ in (\ref{hband}) and $f_0^{(h)}$ in (\ref{optfh}) change places.
\end{remark}

\subsection{Improved reconstruction of nonlinear modes}\label{sec:inmrec}

Given the reconstructed amplitudes $A^{(h)}(t)$, phases $\phi^{(h)}(t)$ and frequencies $\nu^{(h)}(t)$ of all the true harmonics, the most straightforward way to reconstruct the full NM is just to add all $A^{(h)}\cos\phi^{(h)}(t)$ together. However, in this way the NM picks up noise contributions from all harmonics, which can make it quite inaccurate.

Fortunately, there is a cleverer way to perform the reconstruction, also yielding more accurate parameters for the individual harmonics. Thus, one can utilize the theoretical amplitude, phase and frequency relationships between the harmonics, i.e.\ $A^{(h)}=a_hA^{(1)}(t)$, $\phi^{(h)}-h\phi^{(1)}=\varphi_h$ and $\nu^{(h)}(t)=h\nu^{(1)}(t)$. Then, because the components with the higher amplitudes are expected to be less noise-corrupted, one can refine the parameters of each harmonic by weighted averaging over the parameters of all harmonics:
\begin{equation}\label{recNM}
\begin{aligned}
&\tilde{A}^{(h)}(t)=\langle A^{(h)}(t)\rangle\frac{\sum_{h'}A^{(h')}(t)}{\sum_{h'}\langle A^{(h')}(t)\rangle}\\
&\tilde{\phi}^{(h)}(t)=\arg\Big[\sum_{h'}\min(1,h'/h)\langle A^{(h')}(t)\rangle
e^{i(h\phi^{(h')}(t)-\Delta\phi_{h',h}-2\pi I[(h\phi^{(h')}(t)-h'\phi^{(h)}(t)-\Delta\phi_{h',h})/2\pi])/h'}\Big],\\
&\tilde{\nu}^{(h)}(t)=\frac{\sum_{h'}\min(1,h'/h)\langle A^{(h')}(t)\rangle h\nu^{(h')}/h'}{\sum_{h'}\min(1,h'/h)\langle A^{(h')}(t)\rangle}\\
\end{aligned}
\end{equation}
where $\Delta\phi_{h',h}\equiv\arg\big[e^{i(h\phi^{(h')}(t)-h'\phi^{(h)}(t))}\big]\in[-\pi,\pi]$, and $I[\dots]$ denotes rounding to the nearest integer, so that $I[0.8]=1$, $I[-0.6]=-1$ (the corresponding term is needed to eliminate possible noise-induced phase slips, i.e.\ the rapid growth by $2\pi$ in the phase differences). Note also the multiplier $\min(1,h'/h)$, appearing for phase and frequency refinement in (\ref{recNM}). It is needed to account for the scaling of phase and frequency errors of lower harmonics when they are mapped to higher ones. Thus, if $\nu^{(1)}(t)$ has an error $\varepsilon(t)$, then $h\nu^{(1)}(t)$ will have error $h\varepsilon(t)$.

\begin{remark}
The refinement (\ref{recNM}) assumes that the reconstruction error for each harmonic is directly proportional to its amplitude. However, this is true only if, firstly, there are no side components with which harmonics interfere and, secondly, the amount of noise picked up while reconstructing harmonics is the same for each of them. These criteria are rarely both satisfied in practice, so in general one should take into account the mean absolute reconstruction errors $c(h)$ of the $h^{\rm th}$ harmonic parameters. This can be done by changing $A^{(h')}(t)\rightarrow c^{-1}(h')A^{(h')}(t)$ in all the expressions (\ref{recNM}). Unfortunately, the errors $c(h)$ are very hard to estimate, even roughly, and therefore we do not utilize them.
\end{remark}

The refinement (\ref{recNM}) not only makes the reconstructed NM much more accurate, solving the problem of picking up the cumulative noise of all the harmonics, but also reduces the noise level in each harmonic separately. Thus, noise-induced variations in the parameters of different harmonics are expected to be mutually independent, and so they average out being added together. As a result, the more harmonics are contained in the NM, the more accurately it can be reconstructed. While NMD is generally noise-robust, due to being based on time-frequency methods for which only the spectral power of noise in the vicinity of components' frequencies matters, NM reconstruction by (\ref{recNM}) raises this robustness to an extreme extent.

\begin{remark}
It should be noted that the refinement (\ref{recNM}) is not valid if the NM waveform changes in time, i.e.\ the relationships between amplitudes of the harmonics and their phase shifts are non-constant. In this case one should just sum all extracted harmonics, even though some of them might be highly corrupted. Nevertheless, in this work we assume that the waveform is constant, as otherwise even the harmonic identification scheme discussed in Sec.\ \ref{sec:htest} might become inaccurate (if the waveform changes slowly enough, however, it still works well).
\end{remark}

\subsection{Speed improvement}\label{sec:ispeed}

To extract a particular harmonic candidate, one does not need to calculate the TFR over the whole frequency range, which would be computationally expensive. Rather, one can calculate it only in the range where the possible harmonic is expected to lie. As discussed previously, given the frequency range $[\omega_{\min}^{(1)},\omega_{\max}^{(1)}]\equiv[\overline{\omega}_c^{(1)}-\Delta\omega_c^{(1)},\overline{\omega}_c^{(1)}+\Delta\omega_c^{(1)}]$ where the first harmonic is concentrated in TFR, it is clear from (\ref{hband}) that the appropriate frequency range for the $h^{\rm th}$ harmonic will be characterized by $\omega_{\min,\max}^{(h)}=h\overline{\omega}_c^{(1)}\mp\max(1,h)\Delta\omega_c^{(1)}$. The multiplier $\max(1,h)$ takes into account that $h$ can be smaller than unity, since we use $h=1/2,1/3,\dots$ when checking whether the extracted component is first (or higher) harmonic (see Sec.\ \ref{sec:hissues}); the frequency range should not be squeezed for $h<1$ because e.g.\ for NM with only amplitude modulation all harmonics will have the same bandwidth (see Sec.\ \ref{sec:ihadapt}).

Given the extracted time-frequency support of the first harmonic $[\omega_-(t),\omega_+(t)]$, the required frequency range for it to be fully contained in the TFR can be estimated as $[\omega_{\min}^{(1)},\omega_{\max}^{(1)}]=[\min\omega_-(t),\max\omega_+(t)]$. However, the maximum value of the $\omega_+(t)$ might be excessively large (it is in theory infinite e.g.\ if the signal is represented by a single tone); the same applies to minimum $\omega_-(t)$. Therefore, for determination of the harmonic frequency range it is better to use a support $[\tilde{\omega}_-(t),\tilde{\omega}_+(t)]$ that is narrower in frequency, determined as the maximum region where the TFR amplitude of the first harmonic is non-negligible. Mathematically it can be defined as
\begin{equation}\label{etfsupp}
\begin{aligned}
&\tilde{\omega}_{-}(t):\;|H_s(\omega\in[\omega_-(t),\tilde{\omega}_-(t)],t)|<\epsilon_s |H_s(\omega_p(t),t)|\\
&\tilde{\omega}_{+}(t):\;|H_s(\omega\in[\tilde{\omega}_+(t),\omega_+(t)],t)|<\epsilon_s |H_s(\omega_p(t),t)|\\
\end{aligned}
\end{equation}
where $\omega_p(t)\in[\omega_-(t),\omega_+(t)]$ is the position of the peak, as always. If $\tilde{\omega}_-(t)$ (\ref{etfsupp}) does not exist, e.g.\ if the $|H_s(\omega_-(t),t)|\geq\epsilon_s |H_s(\omega_p(t),t)|$, we assign $\tilde{\omega}_-(t)=\omega_-(t)$; we do the same for $\tilde{\omega}_+(t)$ in a similar case. In this work we set $\epsilon_s=0.001$.

Using (\ref{etfsupp}), the frequency range for the $h^{\rm th}$ harmonic can be estimated as
\begin{equation}\label{fharmf}
\begin{gathered}
\omega_{\min,\max}^{(h)}=h\frac{\omega_{\max}^{(1)}+\omega_{\min}^{(1)}}{2}\mp\max(1,h)\frac{\omega_{\max}^{(1)}-\omega_{\min}^{(1)}}{2},\\
\omega_{\min}^{(1)}=\min\big[\underset{0.05}{\rm perc}[\tilde{\omega}_-(t)],{\min}[\omega_p(t)]\big],\;
\omega_{\max}^{(1)}=\max\big[\underset{0.95}{\rm perc}[\tilde{\omega}_+(t)],{\max}[\omega_p(t)]\big],
\end{gathered}
\end{equation}
where ${\rm perc}_{p}$ denotes the $p$'th percentage largest value of the argument (with ${\rm perc}_0$ and ${\rm perc}_1$ corresponding to usual minimum and maximum, respectively). The 95\% largest value of $\tilde{\omega}_+(t)$ instead of the overall maximum is taken in (\ref{fharmf}) to prevent selection of an excessively wide region, because due to noise $\tilde{\omega}_+(t)$ might sometimes be too large; the same applies to $\tilde{\omega}_-(t)$. In the process of harmonic extraction for $h>1$ ($h<1$) it is also useful to restrict $\omega_{\min}^{(h)}$ ($\omega_{\max}^{(h)}$) to be higher than the 5\% (smaller than the 95\%) value of the instantaneous frequency of the last true harmonic being extracted.

The range $[\omega_{\min}^{(h)},\omega_{\max}^{(h)}]$ can, however, be more compressed or stretched, depending on the window/wavelet resolution parameter $f_0^{(h)}$ used for the analysis of the $h^{\rm th}$ harmonic. The expression (\ref{fharmf}) is optimal only for
\begin{equation}\label{resfropt}
\mbox{\textbf{WFT:}}\;f_0^{(h)}=h^{-1}f_0^{(1)}\min(1,h);\;\;\;\mbox{\textbf{WT:}}\;f_0^{(h)}=f_0^{(1)}\min(1,h).
\end{equation}

Indeed, as discussed previously, if the FT of the first harmonic is contained in the band $[\omega_1^{(1)},\omega_2^{(1)}]$, the largest band where $h^{\rm th}$ harmonic can be contained within $\omega_{1,2}^{(h)}=h\frac{\omega_2^{(1)}+\omega_1^{(1)}}{2}\mp\max(1,h)\frac{\omega_2^{(1)}-\omega_1^{(1)}}{2}$. But the range $[\omega_{\min}^{(h)},\omega_{\max}^{(h)}]$ where harmonics lie in the TFR is larger due to the finite frequency resolution, being $[\omega_1^{(h)}-\Delta/f_0,\omega_2^{(h)}+\Delta/f_0]$ for the WFT with a Gaussian window (\ref{winwav}) and $[\omega_1^{(h)}e^{-\Delta/f_0},\omega_2^{(h)}e^{\Delta/f_0}]$ for the WT with a lognormal wavelet (\ref{winwav}), where $\Delta$ is a particular constant. As a result, for $f_0^{(h)}>h^{-1}f_0^{(1)}\min(1,h)$ (WFT) or $f_0^{(h)}>f_0^{(1)}\min(1,h)$ (WT), the range (\ref{fharmf}) is too large; although one can in principle ``squeeze'' it in this case, due to noise the appropriate squeezing cannot always be estimated well, so it is better to leave all as it is. On the other hand, for $f_0^{(h)}<h^{-1}f_0^{(1)}\min(1,h)$ (WFT), or $f_0^{(h)}<f_0^{(1)}\min(1,h)$ (WT), one might need a wider range than (\ref{fharmf}), so it should be stretched. Therefore, if (\ref{resfropt}) is not the case, the values (\ref{fharmf}) should be updated as
\begin{equation}\label{frange}
\begin{aligned}
\mbox{\textbf{WFT:}}\;
&\omega_{\min,\max}^{(h)}=\overline{\omega}^{(h)}+\Delta\omega_{\min,\max}^{(h)}\max(1,f_0^{(1)}\min(1,h^{-1})/f_0^{(h)})\\
\mbox{\textbf{WT:}}\;
&\omega_{\min,\max}^{(h)}=\omega_{\min,\max}^{(h)}\Big(\overline{\omega}^{(h)}/\omega_{\min,\max}^{(h)}\Big)^{1-\max(1,f_0^{(1)}\min(1,h)/f_0^{(h)})},\\
\end{aligned}
\end{equation}
where $\overline{\omega}^{(h)}\equiv\frac{h}{2}\big[\omega_{\min}^{(1)}+\omega_{\max}^{(1)}\big]$ and $\Delta\omega_{\min,\max}^{(h)}\equiv \omega_{\min,\max}^{(h)}-\overline{\omega}^{(h)}$ are calculated from the estimates (\ref{fharmf}). The expressions (\ref{fharmf}) and (\ref{frange}) together give the frequency range in which to calculate the WFT/WT (based on the window/wavelet (\ref{winwav}) with the resolution parameter $f_0^{(h)}$) for the $h^{\rm th}$ harmonic.

To reduce the computational cost of the surrogate test against noise, utilized for the stopping criterion (see Sec.\ \ref{sec:stopcrit}), one can use the same trick as above, calculating TFRs for surrogates in the frequency range where the original component resides. After extracting the component from the original TFR, this range can be estimated by (\ref{fharmf}) with $h=1$.

\section{Parameters and their choice}\label{sec:pchoice}

There are different settings that can be used while applying NMD. However, many of the parameters are either set at well-established values, or can be chosen adaptively, thus removing the ambiguity and making NMD a kind of superadaptive method. The settings and their choice are discussed below.

\subsection{Resolution parameter $f_0$}\label{sec:prespar}

The resolution parameter $f_0$ of the window/wavelet (\ref{winwav}) determines the time and frequency resolution of the WFT/WT, i.e.\ how fast time changes can be reflected and how close in frequency components can be resolved, respectively. Time and frequency resolutions are not independent, being inversely proportional to each other, and an increase of $f_0$ increases the frequency resolution but decreases the time resolution. This inability to locate the component precisely in both time and frequency is a very general issue, which manifests itself in nearly every aspect of signal analysis and represents the time-frequency analogue of the Heisenberg uncertainty principle in quantum mechanics \cite{Folland:97,Cohen:95,Kaiser:94,Mallat:08}.

Although we have discussed how to choose $f_0$ for harmonics given the extracted component, it is still unclear how to choose it at the first step, when we extract the main component. Generally, one needs to select $f_0$ based on a compromise between better reflecting time variability and better resolving components in frequency, so the optimal choice depends on the characteristics of the components contained in the signal. Unfortunately, at present there does not seem to be any universal approach enabling one to choose $f_0$ appropriately for any given signal (see \cite{Iatsenko:13tfr2,Mallat:08} for a discussion of this issue and the effects of different choices), so it remains the only important parameter of NMD that cannot be adaptively selected. Its choice, however, seems to be slightly more universal for the WT, as the latter adjusts its resolution in relation to the frequency band being studied; for it, one usually uses $f_0=1$, setting some standard limits on the allowable relative amplitude/frequency modulation of the components and the frequency distances between them.

Nevertheless, the very possibility of adjusting the time and frequency resolution is a great advantage of TFR-based methods. For example, as discussed previously in Remark \ref{rem:emdres}, the (E)EMD has time and frequency resolution properties similar to WT but, in contrast to the latter, its resolutions are set around some implicit values and cannot be changed. The choice of $f_0$ therefore gives NMD more flexibility in comparison to methods not possessing such a parameter.

\subsection{TFR type: WFT or WT?}\label{sec:ptfrtype}

The difference between the WFT and the WT lies in the type of frequency resolution, which is linear for the former and logarithmic for the latter. Thus, the ability of the WT to resolve components in frequency worsens with increasing frequency, while its ability to reflect time variations improves; in contrast, the WFT does not discriminate between components based on their characteristic frequencies. Preference for one type of TFR over the other therefore depends on the signal structure. The WT is especially suitable if the higher frequency components contained in the signal are more distant in frequency and have higher time variability than those at lower frequencies, which is often the case for real signals; otherwise, the WFT is preferable.

Without some {\it a priori} knowledge, it is hard to select adaptively (i.e.\ automatically) the most appropriate type, especially given the associated problems related to the choice of the resolution parameter, as discussed above. However, after one component has been extracted (even roughly), selection of the optimal resolution type can be made based on its properties. Thus, if the time-variability of the component's parameters increases with frequency, then the WT is the most suitable, whereas otherwise one should prefer the WFT. Given the initially extracted component's amplitude $A(t)$, phase $\phi(t)$ and frequency $\nu(t)$, an empirical condition for selecting the most appropriate TFR type can be stated as:
\begin{equation}\label{ctype}
\Big(1+V[\partial_t\nu(t),\nu(t)]\Big)^{-1}+\Big(1+V[\partial_tA(t),\nu(t)]\Big)^{-1}<1\Rightarrow\mbox{use WFT};
\;\;\mbox{otherwise}\Rightarrow\mbox{ use WT},\\
\end{equation}
where
\begin{equation}\label{ctype}
V[x(t),y(t)]\equiv \frac{\std\big[|x(t)/y(t)|^+\big]}{\std\big[|x(t)/\langle y(t)\rangle|^+\big]},\;\;\Big(\;|f(t)|^+\equiv|f^+(t)|\;\Big).\\
\end{equation}
Thus, the values of $V[\partial_tA(t),\nu(t)]$ and $V[\partial_t\nu(t),\nu(t)]$ quantify whether the amplitude and frequency modulation of the component become stronger with increasing frequency ($V<1$) or not ($V>1$). In the former case a reliable representation of the component requires higher time-resolution for higher frequencies, so the WT is preferred, while for the latter case the WFT should be used. For example, for the linear (hyperbolic) chirp $\nu(t)=\nu_0+a t$ ($\nu(t)=\exp(at)$) one has $V[\partial_t\nu(t),\nu(t)]=\infty$ ($V[\partial_t\nu(t),\nu(t)]=0$), reflecting the fact that the WFT (WT) is most appropriate for its representation, in accordance with what is known \cite{Mallat:08}.

The derivatives $\partial_t\nu(t)$ and $\partial_t A(t)$ in (\ref{ctype}) can be estimated by numerical differentiation. However, when noise is present in $A(t)$ and $\nu(t)$, it will be greatly amplified in such estimates, so they will be generally quite noisy. Consequently, instead of standard deviations $\std[x(t)]$ in (\ref{ctype}), it is better to use 75-percentiles, i.e.\ the width of the range within which 75\% of the values of $x(t)$ are contained. Alternatively, one can reconstruct $\partial_t\nu(t)$ and $\partial_t A(t)$ by deriving the direct reconstruction formulas for them as explained in \cite{Iatsenko:13tfr1}.

The remained question is which TFR type to use for the preliminary signal exploration, i.e.\ extraction of the component's parameters to be used in (\ref{ctype}). As was discussed, the answer depends on the signal structure, and there is no universal criterion. However, since the WT is usually faster to calculate (due to its logarithmic scale) and generally has a more universal choice of the resolution parameter, we use it.

Summarizing, we calculate the WT of the signal and extract from it the component and its parameters. Then we utilize the criterion (\ref{ctype}) and determine which is the best type of TFR to use in the given case. If this is the WT, we retain the component already extracted; otherwise, we calculate the WFT of the signal in the corresponding frequency range (\ref{fharmf}) and re-extract all the parameters from it. To preserve the time and frequency resolution properties for which the component was extracted from the WT (which we assume to be appropriate for it), the resolution parameter for the WFT $f_0^{(WFT)}$ should be adjusted accordingly. If the WFT and WT are calculated using a Gaussian window and lognormal wavelet (\ref{winwav}), the rule is
\begin{equation}\label{f0wftwt}
f_0^{(WFT)}\approx 2\pi f_0^{(WT)}/\langle\nu(t)\rangle.
\end{equation}
It follows from the fact that, for window and wavelet (\ref{winwav}) and not too small $f_0^{(WT)}$, the frequency resolution of the WT around $\omega=2\pi$ (if ``linearized'') is similar to the frequency resolution of the WFT for the same resolution parameter \cite{Iatsenko:13tfr1}; taking into account the scaling nature of the WT, i.e.\ its frequency-dependent resolution, one then obtains (\ref{f0wftwt}). Since for all harmonics by definition one would have identical $V[\partial_t\nu(t),\nu(t)]$ and $V[\partial_tA(t),\nu(t)]$ in (\ref{ctype}), the same type of frequency resolution (linear or logarithmic) is appropriate for all of them. Hence, for the extraction of harmonics we use the same TFR type as was selected for the original component.

\begin{remark}
The logarithmic frequency resolution of the WT might by itself introduce correlation between the frequency and the amplitude/frequency variations of the component. For example, when the signal is corrupted with white noise, its power in the WT will be proportional to frequency, leading to larger noise-induced amplitude/frequency variations in the extracted component at higher frequencies. Such artificial correlation might cause the criterion (\ref{ctype}) to select the WT as the preferred representation even when this is not the case. To avoid this, we use threshold 1.1 instead of 1 in (\ref{ctype}).
\end{remark}

\subsection{Reconstruction method: direct or ridge?}\label{sec:precm}

As mentioned in Sec.\ \ref{sec:cextract}, one can use either of two alternative methods for reconstruction of the components from the WFT/WT: direct (\ref{dirrec}) or ridge (\ref{ridgerec}). The differences and errors of each method were studied in detail in \cite{Iatsenko:13tfr2}. It was found that the ridge method is more noise-robust, but that the direct method allows the time variability in the component's parameters at low noise levels to be followed more accurately.

For some parts of the NMD procedure it is inherently better to use a particular reconstruction method. Thus, in the criterion for selecting the TFR type (\ref{ctype}) we always use the parameters reconstructed by the ridge method due to the noise-robustness of such estimates. Additionally, while testing the signal against noise in the stopping criterion (Sec.\ \ref{sec:stopcrit}) we also use the ridge estimates for calculating the discriminating statistics $D$ (\ref{sentropy}). This seems natural because the curve extraction is based on the amplitudes and frequencies of the peaks (see Sec.\ \ref{sec:cextract}), and we have found such an approach to be more stable than using the direct estimates; the noise-robustness of ridge reconstruction is advantageous here as well.

For the other parts of the procedure, the reconstruction method can be chosen adaptively. To choose it at the first step, when we have extracted the time-frequency support $[\omega_-(t),\omega_+(t)]$ of some component and need to reconstruct its parameters, one can use the approach suggested in \cite{Iatsenko:13tfr2}. Thus, we first reconstruct the parameters by both methods, getting $A^{(d,r)}(t),\phi^{(d,r)}(t),\nu^{(d,r)}(t)$, where ``d'' and ``r'' denote direct and ridge estimates, respectively. Then we compute the TFR (using the same parameters as originally) of the signal $s^{(d)}(t)=A^{(d)}\cos\phi^{(d)}$, extract the ridge curve from it (taking simple maximuma $\omega_p(t)=\argmax_\omega|H_s(\omega,t)|$ is sufficient here), and reconstruct by direct method the ``refined'' parameters $\tilde{A}^{(d)},\tilde{\phi}^{(d)},\tilde{\nu}^{(d)}$. We do the same for the ``ridge'' signal $s^{(r)}(t)=A^{(r)}\cos\phi^{(r)}$, now using ridge method to reconstruct the refined estimates. Then we compute inconsistencies $\varepsilon_{a,\phi,\nu}^{(d,r)}$ between the original and refined estimates for each method:
\begin{equation}\label{apfrel}
\begin{aligned}
&\varepsilon_A^{(d,r)}\equiv \kappa_A^{(d,r)}\sqrt{\langle (\tilde{A}^{(d,r)}(t)-A^{(d,r)}(t))^2 \rangle},\\
&\varepsilon_\phi^{(d,r)}\equiv \kappa_\phi^{(d,r)}\sqrt{1-|\langle e^{i(\tilde{\phi}^{(d,r)}(t)-\tilde{\phi}^{(d,r)}(t))} \rangle|^2},\\
&\varepsilon_\nu^{(d,r)}\equiv \kappa_\nu^{(d,r)}\sqrt{\langle (\tilde{\nu}^{(d,r)}(t)-\nu^{(d,r)}(t))^2\rangle},\\
\end{aligned}
\end{equation}
where $\kappa_{A,\phi,\nu}^{(d,r)}$ are the coefficients to tune the performance of the approach (we empirically found $\kappa_{A,\phi,\nu}^{(d)}=\{3,4,2\}$, $\kappa_{A,\phi,\nu}^{(r)}=1$). Obviously, for the exact estimates one has $\varepsilon_{A,\phi,\nu}=0$, so it is natural to assess the relative performance of the reconstruction methods based on their inconsistencies (\ref{apfrel}). Therefore, we use the direct amplitude estimate if $\varepsilon_A^{(d)}<\varepsilon_A^{(r)}$, and the ridge amplitude estimate otherwise; the same applies to the phase and frequency estimates. Although empirical, this approach works very well in practice.

\begin{remark}
To implement the criterion (\ref{apfrel}), one needs to compute two new TFRs, for the direct and ridge signals $s^{(d,r)}(t)$. Obviously, one does not need to calculate them for all frequencies, and one can restrict the frequency range to that of the original component $[\omega_{\min}^{(h)},\omega_{\max}^{(h)}]$, which can be estimated by (\ref{fharmf}). This will improve the speed of the procedure.
\end{remark}

For harmonics the appropriate reconstruction method can be determined simply as that giving the highest consistency $\rho^{(h)}$ (\ref{apc}). Therefore, while adapting the resolution parameter for the harmonic, at each $f_0$ we reconstruct its parameters by both methods and at the end select the one characterized by highest $\rho^{(h)}$.

An adaptive choice of reconstruction method, as outlined above, ensures that in most cases we get the best possible estimates: when the amplitude/frequency modulation is low and noise is high the noise-robust ridge reconstruction method is used, while in the opposite case the direct method is applied, giving more accurate amplitude/frequency variations. This approach further increases the noise-robustness and the accuracy of the NMD.

\subsection{Other parameters}\label{sec:pother}

All the other parameters can be partitioned onto two groups: some pre-fixed settings and the parameters of numerical accuracy. The former group includes the significance levels for the two surrogate tests (used for identification of the harmonics, Sec.\ \ref{sec:htest}, and for the stopping criterion, Sec.\ \ref{sec:stopcrit}), the minimal consistency $\rho_{\min}$ (\ref{rhomin}) etc. Each of these parameters is set to a well-established value, either corresponding to a standard convention (such as $95\%$ significance of the surrogate tests \cite{Schreiber:00,Theiler:92}) or found empirically (e.g.\ the expression (\ref{rhomin}) for $\rho_{\min}$).

The second group includes the accuracy with which to determine the optimal $f_0$ while adapting it to the harmonics, the precision $\epsilon_s$ which to use for determining the minimal support (\ref{etfsupp}) etc. Here one faces the usual tradeoff between accuracy and computational cost. The default values, however, are sufficient in most cases and further increase of precision leads to only a slight improvement in the accuracy of the method.

\section{Simulation examples}\label{sec:simex}

Having described all parts of the NMD procedure, we now illustrate the method by consideration of some specific examples. In this section we consider simulated signals whose composition is precisely known, so that one can assess reliably the performance of NMD and compare it to that of the other methods.

\subsection{Example 1}\label{sec:simex1}

\begin{figure}[t]
\includegraphics[width=1.0\linewidth]{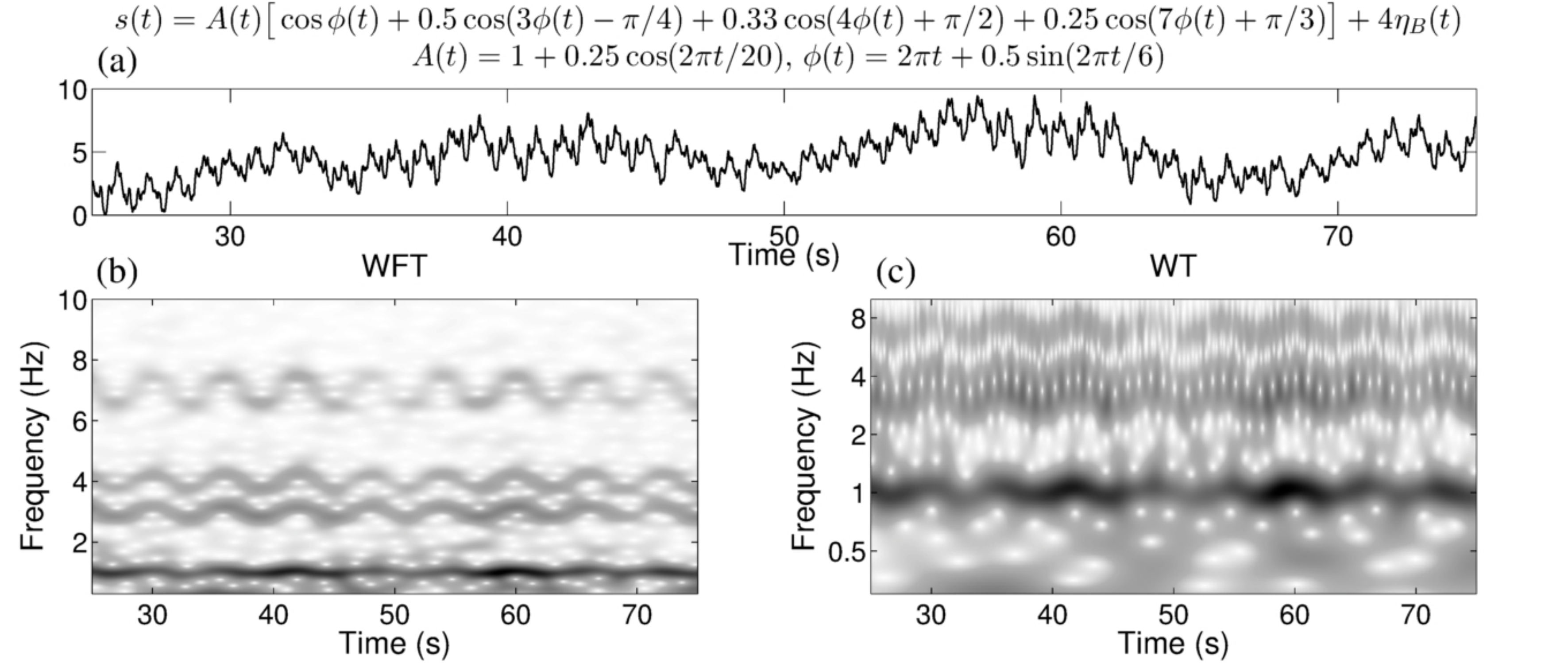}
\caption{(a) The central 50\,s section of the signal $s(t)$ specified by the equation at the top, which represents a single NM corrupted by Brownian noise of standard deviation equal to $4$ (Brownian noise of unit deviation $\eta_B(t_n)\sim \sum_{m=1}^n\eta_W(t_m)$ is obtained as the normalized cumulative sum of the Gaussian white noise signal $\eta_W(t)$). (b),(c) Central parts of the WFT and WT, respectively, of the signal $s(t)$ shown in (a). The signal was sampled at 100 Hz for 100 s.}\label{fig:ex1tfr}
\end{figure}

As a first, simple and illustrative example, we take the signal depicted in Fig.\ \ref{fig:ex1tfr}(a). Its WFT and WT are shown in Fig.\ \ref{fig:ex1tfr}(b) and (c), respectively. From the figure it is immediately evident that, while in the WFT the harmonics with frequencies around 3 and 4 Hz are distinguishable, in the WT they interfere strongly and cannot be resolved; in contrast, the WT has high enough time-resolution at 7 Hz to represent the frequency modulation of the corresponding harmonic, while in the WFT this highest harmonic self-interferes (there appear the ``bridges'' between consecutive frequency modulation cycles), indicating that the time resolution is insufficient. Therefore, for the present case one cannot appropriately extract all harmonics using either WFT or WT with constant $f_0$. However, adaptive representation of the harmonics, as discussed in Sec.\ \ref{sec:ihadapt}, solves the problem.

The NMD proceeds as follows. First it extracts and reconstructs the dominant component from the signal's WT (located at 1 Hz and corresponding to the first harmonic of the NM in the present case). It is tested with the surrogates against noise (see Sec.\ \ref{sec:stopcrit}) and passes the test (100\% significance). Using (\ref{ctype}), the WFT is determined to be more suitable for the representation of this component, because its amplitude and frequency modulations do not depend on frequency; this type of TFR is then used in what follows. The component is re-extracted from the signal's WFT calculated in the corresponding frequency range (around 1 Hz), and its parameters are reconstructed by both direct (\ref{dirrec}) and ridge (\ref{ridgerec}) methods. Using (\ref{apfrel}), it is established that ridge estimates of amplitude, phase and frequency seem to be more accurate in the present case. This is to be expected, given the high noise level around 1 Hz and the pronounced noise-susceptibility of direct estimates. The ridge estimates are therefore taken as the ones to be used.

Next, we test whether the extracted component is the first harmonic by checking its possible $1/h$ harmonics. Thus, first the $1/2$ harmonic is extracted and tested for being true (see Sec.\ \ref{sec:htest}) using different resolution parameters $f_0$ within the range (\ref{optfh}); among the direct and ridge estimates obtained for $f_0$ at which the harmonic was identified as true, we choose those maximizing the consistency (\ref{apc}). In the present case, the $1/2$ harmonic is identified as false for all tested $f_0$, so it is discarded. We do the same for 1/3 and 1/4 harmonics, which are both identified as false. Since there are $S_0=3$ consecutive false harmonics (1/2, 1/3 and 1/4), we stop the procedure and correctly conclude that the extracted component is first (and not higher) harmonic of the NM.

\begin{figure}[t]
\includegraphics[width=1.0\linewidth]{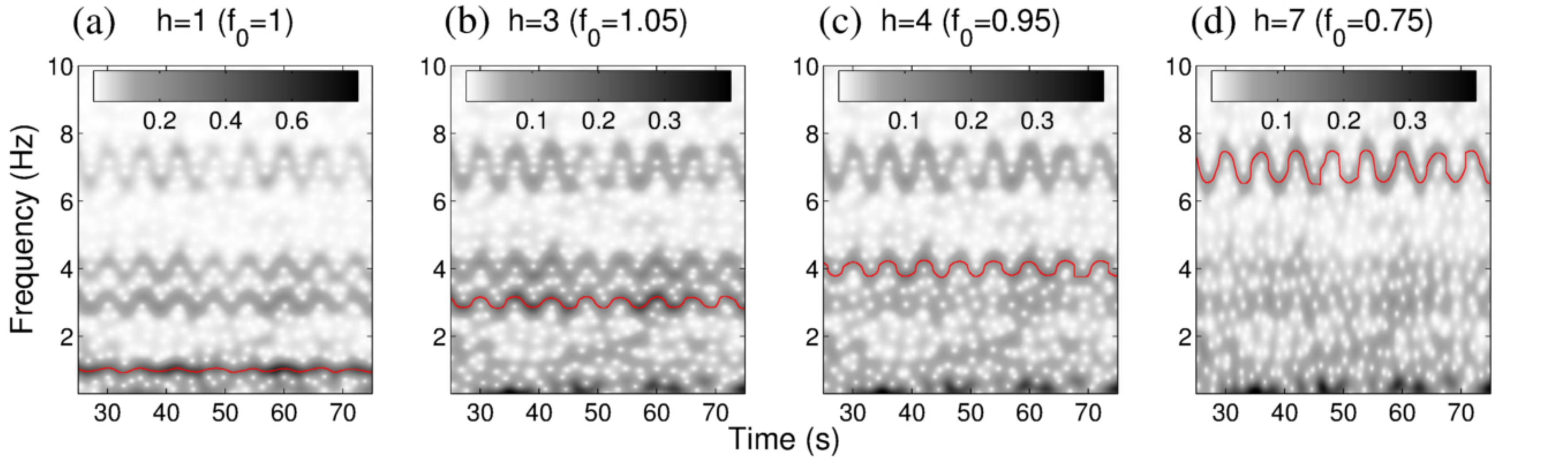}
\caption{The WFTs from which each true harmonic is reconstructed, with the corresponding extracted ridge curves being shown by solid red lines. The window resolution parameter $f_0$ is adjusted individually for each harmonic. After the harmonic is reconstructed, it is subtracted from the signal, so that e.g.\ the first harmonic no longer appears in the WFTs of (b)-(d). Note that the color scaling in (b)-(d) covers half the amplitude range of that in (a).}
\label{fig:ex1harm}
\end{figure}

We then extract and test the higher harmonics $h=2,3,...$ in qualitatively the same way as was done for $h=1/2,1/3,...$. If some harmonic is identified as true, we subtract it from the signal to remove its interference with higher harmonics. As a result of the procedure, all genuine harmonics $h=3,4,7$ were correctly identified as true, and all the others as false. The resolution parameters were adapted for each harmonic so as to optimally represent and reconstruct it, as discussed in Sec.\ \ref{sec:ihadapt} and illustrated in Fig.\ \ref{fig:ex1harm} for the present case. Harmonic extraction was stopped when $S=3$ consecutive harmonics $h=8,9,10$ were identified as false.

\begin{figure}[t!]
\includegraphics[width=1.0\linewidth]{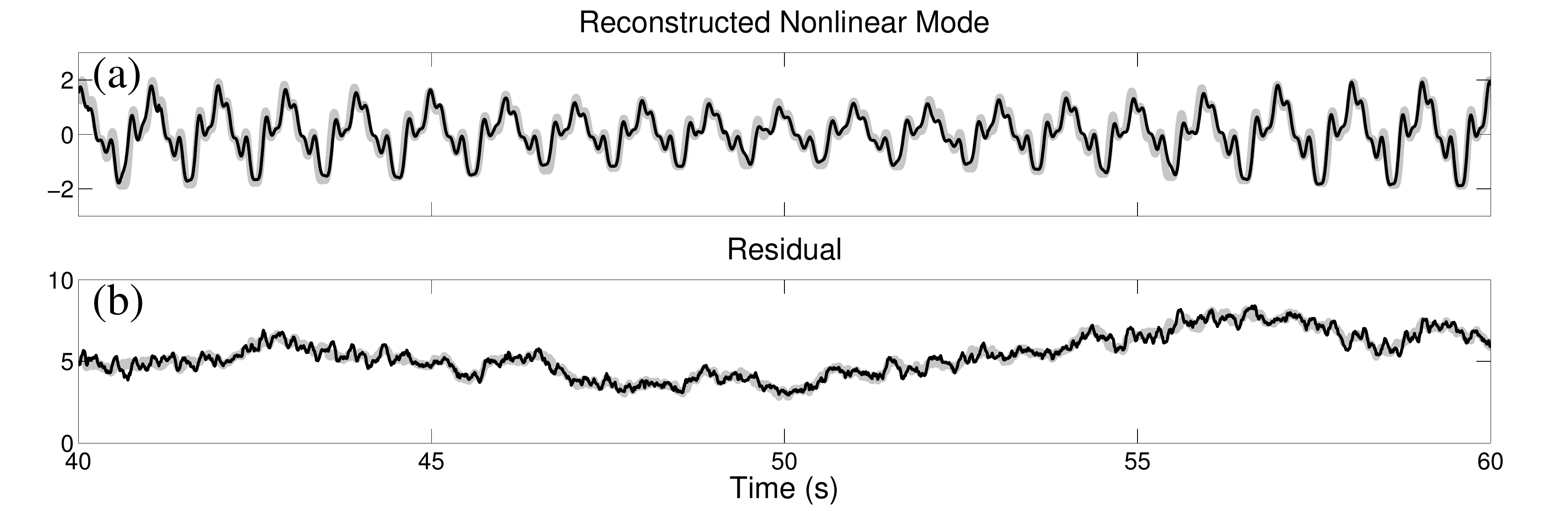}
\caption{The result of applying NMD to the signal shown in Fig.\ \ref{fig:ex1tfr}(a): (a) the reconstructed NM (black line) is compared with the true NM (gray background line); (b) the residual provided by NMD (black line) is compared with the actual background noise (gray background line).}
\label{fig:ex1nmd}
\end{figure}

We then refine the harmonics' parameters using (\ref{recNM}) and reconstruct the full Nonlinear Mode. This NM is then subtracted from the signal, and the procedure is repeated. However, when we extract the next component, it does not pass the surrogate test against noise (see Sec.\ \ref{sec:stopcrit}), and therefore NMD is stopped, with one NM being extracted and the residual correctly identified as noise (in the present case Brownian). The result of NMD is shown in Fig.\ \ref{fig:ex1nmd}, from which one can see that even the residual Brownian noise is well-recovered. To the best of the authors' knowledge, there is at present no other method that can decompose even the current relatively simple signal (shown in Fig.\ \ref{fig:ex1tfr}(a)) in such an accurate and physically meaningfully way.

\begin{figure}[t!]
\includegraphics[width=1.0\linewidth]{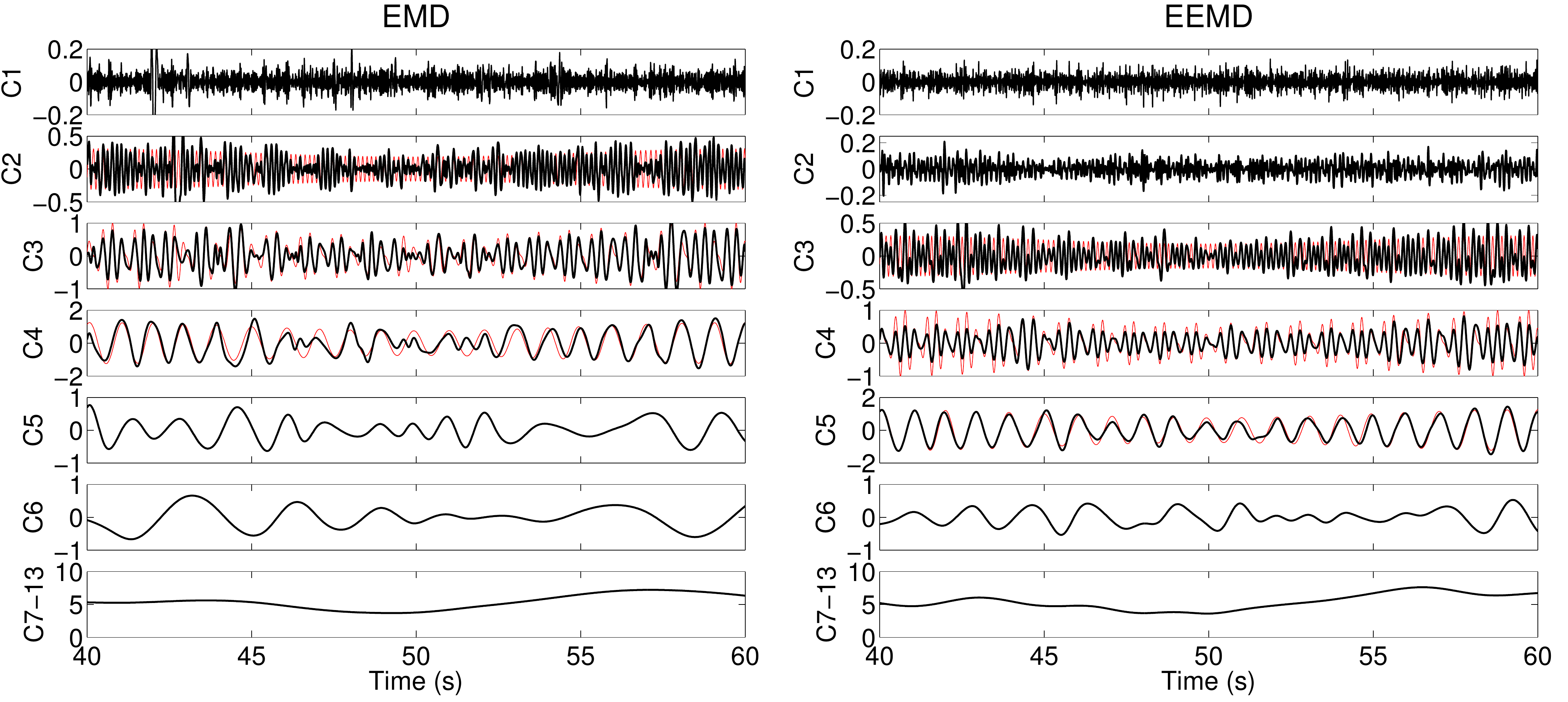}
\caption{The result of applying EMD (left panels) and EEMD (right panels) to the signal shown in Fig.\ \ref{fig:ex1tfr}. Red lines, where present, show the real 1st harmonic (in C4 for EMD and C5 for EEMD), sum of the 3rd and 4th harmonics (in C3 for EMD and C4 for EEMD) and the 7th harmonic (in C2 for EMD and C3 for EEMD). The bottom panels show the sum of the components 7 to 13. For EEMD, we have used 1000 noise realizations with standard deviations ten times smaller than that of the signal.}
\label{fig:ex1emd}
\end{figure}

For example, the results of EMD \cite{Huang:98} and EEMD \cite{Wu:09} procedures are shown in Fig.\ \ref{fig:ex1emd}. In contrast to NMD, (E)EMD produces 13 distinct components, with only the first harmonic of the NM being more-or-less reliably recovered in one of them. Thus, in Fig.\ \ref{fig:ex1emd} C4 for EMD and C5 for EEMD represent the first harmonic, C3 for EMD and C4 for EEMD is the noise-spoiled mix of the 3 and 4 harmonics, C2 for EMD and C3 for EEMD is the badly corrupted 7th harmonic (with influence from harmonics 3 and 4 in the case of EEMD), while none of the other ``components'' has any physical meaning at all.

\subsection{Example 2}\label{sec:simex2}

\begin{figure}[t!]
\includegraphics[width=1.0\linewidth]{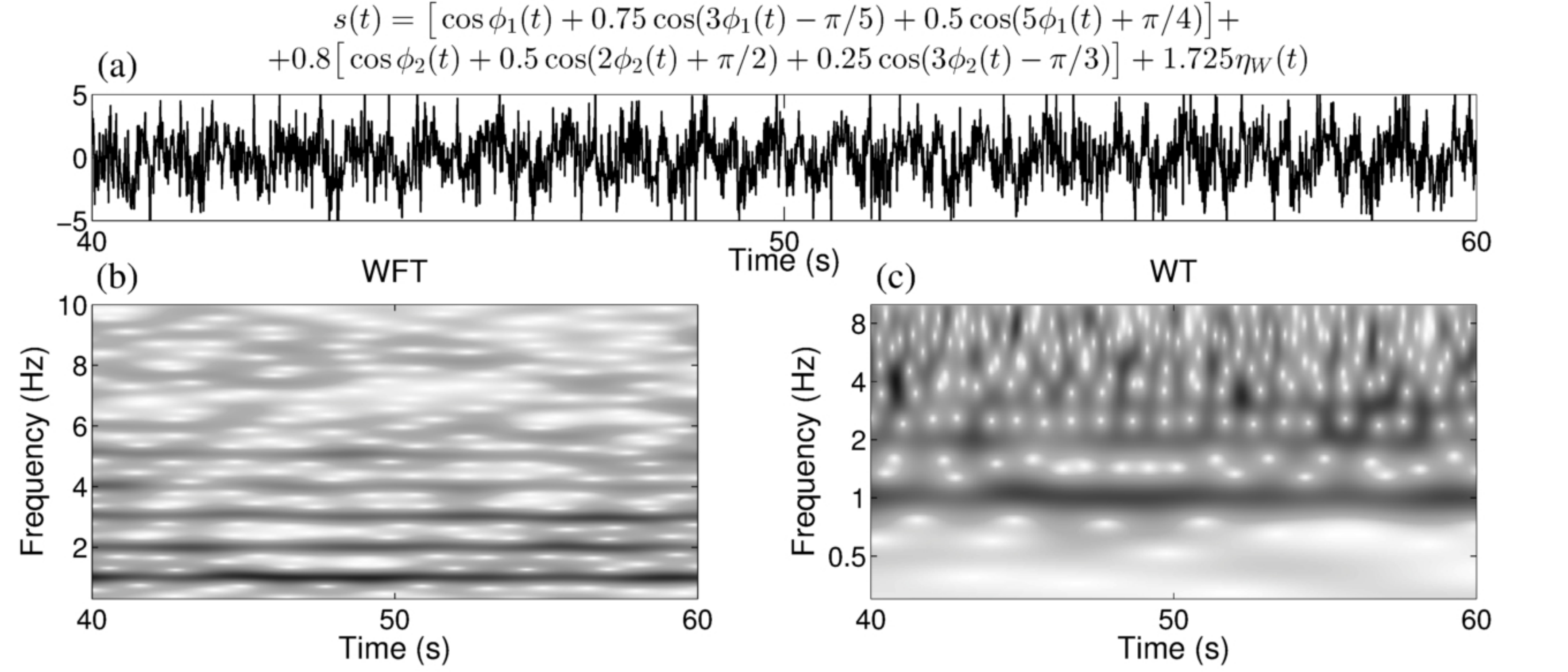}
\caption{(a) The central 20\,s section of the signal $s(t)$ specified by the equation at the top, representing the sum of two NMs corrupted by Gaussian white noise of standard deviation equal to $1.725$; the phases of the NMs $\phi_{1,2}(t)$ were obtained as $\phi_{1,2}(t)=\int_0^t \nu_{1,2}(\tau)d\tau$ with $\nu_{1,2}(t)/2\pi=(1,2)+0.01\tilde{\eta}_{B;1,2}(t)$, where $\tilde{\eta}_{B;1,2}(t)$ are two independent realizations of the unit deviation Brownian noise filtered in the range $[0.01,0.2]$ Hz. (b),(c) Central parts of the WFT and WT of the signal $s(t)$ shown in (a). The signal was sampled at 100 Hz for 100 s.}\label{fig:ex2tfr}
\end{figure}

To demonstrate the exceptional noise-robustness of NMD and the power of the surrogate test in distinguishing true from false harmonics, we consider the signal shown together with its WFT and WT in Fig.\ \ref{fig:ex2tfr}. As can be seen, the harmonics of the second NM are located exactly at the places where the harmonics of the first NM are expected to be, so that they can very easily be confused with them. Moreover, it is very hard to distinguish true from false harmonics in the present case because each NM has constant amplitudes and only very small frequency modulation (the absolute deviation between the expected frequency of the second harmonic of the first NM and the frequency of the first harmonic of the second NM is only $|2\nu_1^{(1)}(t)-\nu_2^{(1)}|/2\pi=0.016\pm0.014$ Hz). Furthermore, the noise is exceedingly strong, in standard deviation being 1.5 times that of the full noise-free signal, 1.8 times that of the first NM, 2.7 times that of the second NM, and 12.2 times that of the smallest (3rd) harmonic of the second NM, located at around 6 Hz (the latter is buried under the noise and not even visible in the WFT shown in Fig.\ \ref{fig:ex2tfr}(b)).

Because the noise is white in the present case, its power in the WT grows with frequency, as seen from Fig.\ \ref{fig:ex2tfr} (note, however, that such a situation rarely arises for real signals). Consequently, instead of using the WT and then adaptively selecting the appropriate type of TFR, for the present case we use the WFT from the very beginning. Nonlinear Mode Decomposition then proceeds as usual: it first extracts the dominant component and tests it against noise; if it passes the test, NMD extracts its harmonics and identifies the true ones; then the full NM is reconstructed and subtracted from the signal, after which the procedure is repeated on the residual until it does not pass the surrogate test against noise.

The relevant information about the NMD outcome is summarized in Table \ref{tab:harmonics}. As can be seen, NMD correctly identifies all harmonics and joins them into two distinct modes. Moreover, the amplitude ratios $a_h$ and phase shifts $\varphi_h$ for each NM (see (\ref{nm})), calculated from the reconstructed amplitude and phases as $a_h=\langle A^{(h)}\rangle/\langle A^{(1)}\rangle$ and $\varphi_h\equiv\arg\big[\langle e^{i(\phi^{(h)}-h\phi^{(1)})}\rangle\big]$, are very close to their actual values; this is true even for the (buried in noise) 3rd harmonic of the second NM. The ridge method is automatically selected for reconstructing of all harmonics, which is indeed the preferable choice due to its noise-robustness \cite{Iatsenko:13tfr2}.

From Table 1 it can be noted that, for all true harmonics, the resolution parameter $f_0$ used is higher than the original. This is because the higher the $f_0$, the easier it is to segregate the component from noise \cite{Iatsenko:13tfr2}. However, increasing $f_0$ at the same time worsens the accuracy of representation of amplitude/frequency modulation of the component \cite{Iatsenko:13tfr2}, so its choice is determined by a compromise between reflecting well the time-variability of component's parameters and suppressing the noise; the adaptive scheme that we use, described in Sec.\ \ref{sec:ihadapt}, effectively implements this criterion.

\begin{center}
\begin{table}[t!]
\small
  \begin{tabular}{| >{\centering\arraybackslash}p{1.25cm} || >{\centering\arraybackslash}p{1cm} | >{\centering\arraybackslash}p{1.2cm} | >{\centering\arraybackslash}p{1cm} | >{\centering\arraybackslash}p{1cm} | >{\centering\arraybackslash}p{1cm} | >{\centering\arraybackslash}p{1cm} | >{\centering\arraybackslash}p{1cm} | >{\centering\arraybackslash}p{1cm} |}
  \hline
  \multirow{2}{1.25cm}{Harmonic number $h$} & \multirow{2}{*}{$\rho^{(h)}$} & signif. & rec. & \multirow{2}{*}{$f_0^{(h)}$} & \multicolumn{2}{c|}{amp. ratio $a_h$} & \multicolumn{2}{c|}{phase shift $\varphi_h/\pi$} \\ \cline{6-9}
                             &               &     level     &      method     &                &    true     &    extr.    &    true     &    extr     \\
                             \hline
  \multicolumn{9}{|c|}{First NM (significance against noise = 100\%)} \\ \hline
  \rowcolor{shadegray}
  \textbf{1}                 &      ---      &      ---      &      ridge      &      1.00      &      1      &      1      &      0      &      0      \\
  2                          &     0.41      &      50\%     &      ridge      &      2.85      &     ---     &     ---     &     ---     &     ---     \\
  \rowcolor{shadegray}
  \textbf{3}                 &     0.75      &     100\%     &      ridge      &      2.08      &     0.75    &     0.74    &    -0.20    &    -0.20    \\
  4                          &     0.29      &      69\%     &      ridge      &      1.77      &     ---     &     ---     &     ---     &     ---     \\
  \rowcolor{shadegray}
  \textbf{5}                 &     0.63      &      98\% &      ridge      &      2.01      &     0.50    &     0.51    &     0.25    &     0.26    \\
  6                          &     0.11      &      52\%     &      ridge      &      1.77      &     ---     &     ---     &     ---     &     ---     \\
  7                          &     0.04      &      67\%     &      direct     &      1.21      &     ---     &     ---     &     ---     &     ---     \\
  8                          &     0.06      &      67\%     &      ridge      &      0.41      &     ---     &     ---     &     ---     &     ---     \\
  \hline
  \multicolumn{9}{|c|}{Second NM (significance against noise = 100\%)} \\ \hline
  \rowcolor{shadegray}
  \textbf{1}                 &      ---      &      ---      &      ridge      &      1.00      &      1      &      1      &      0      &      0      \\
  \rowcolor{shadegray}
  \textbf{2}                 &      0.68     &      98\% &      ridge      &      2.85      &     0.50    &     0.45    &     0.50    &     0.51    \\
  \rowcolor{shadegray}
  \textbf{3}                 &      0.40     &      95\%     &      ridge      &      1.71      &     0.25    &     0.29    &    -0.33    &    -0.33    \\
  4                          &      0.09     &      79\%     &      ridge      &      1.95      &     ---     &     ---     &     ---     &     ---     \\
  5                          &      0.11     &      91\%     &      ridge      &      0.64      &     ---     &     ---     &     ---     &     ---     \\
  6                          &      0.07     &      17\%     &      ridge      &      1.00      &     ---     &     ---     &     ---     &     ---     \\
  \hline
  \multicolumn{9}{|c|}{Residual (significance against noise = 37\%)} \\ \hline
  \end{tabular}
\caption{Summary of the results of NMD applied to the signal shown in Fig.\ \ref{fig:ex2tfr}. For each NM the significance of the surrogate test against noise (see Sec.\ \ref{sec:stopcrit}) was based on 100 surrogates. The columns left-to-right provide information for the $h^{\rm th}$ harmonic: the value of the consistency $\rho^{(h)}(1,1,0)$ (\ref{apc}); the significance of the test against independence (see Sec.\ \ref{sec:htest}); the method chosen for reconstruction of the harmonic (see Sec.\ \ref{sec:precm}); the resolution parameter $f_0^{(h)}$ adapted for the harmonic considered (see Sec.\ \ref{sec:ihadapt}); the true amplitude ratio $a_h\equiv A^{(h)}(t)/A^{(1)}(t)$; the extracted amplitude ratio; the true phase shift $\varphi_h\equiv\phi^{(h)}(t)-h\phi^{(1)}(t)$; the extracted phase shift. Harmonic was identified as true if both $\rho^{(h)}\geq\rho_{\min}=0.25$ (\ref{rhomin}) and the significance level is $\geq95\%$, with the corresponding rows being shaded with gray. The numbers of real harmonics are emboldened, and it can be seen that only they were identified as being true. For simplicity, the results for $h=1/n$ harmonics (which are tested to check whether the extracted component is the first harmonic) are not shown; in each case all three consecutive $h=1/2,1/3,1/4$ harmonics were identified as false.}
\label{tab:harmonics}
\end{table}
\end{center}

\begin{figure}[t]
\includegraphics[width=1.0\linewidth]{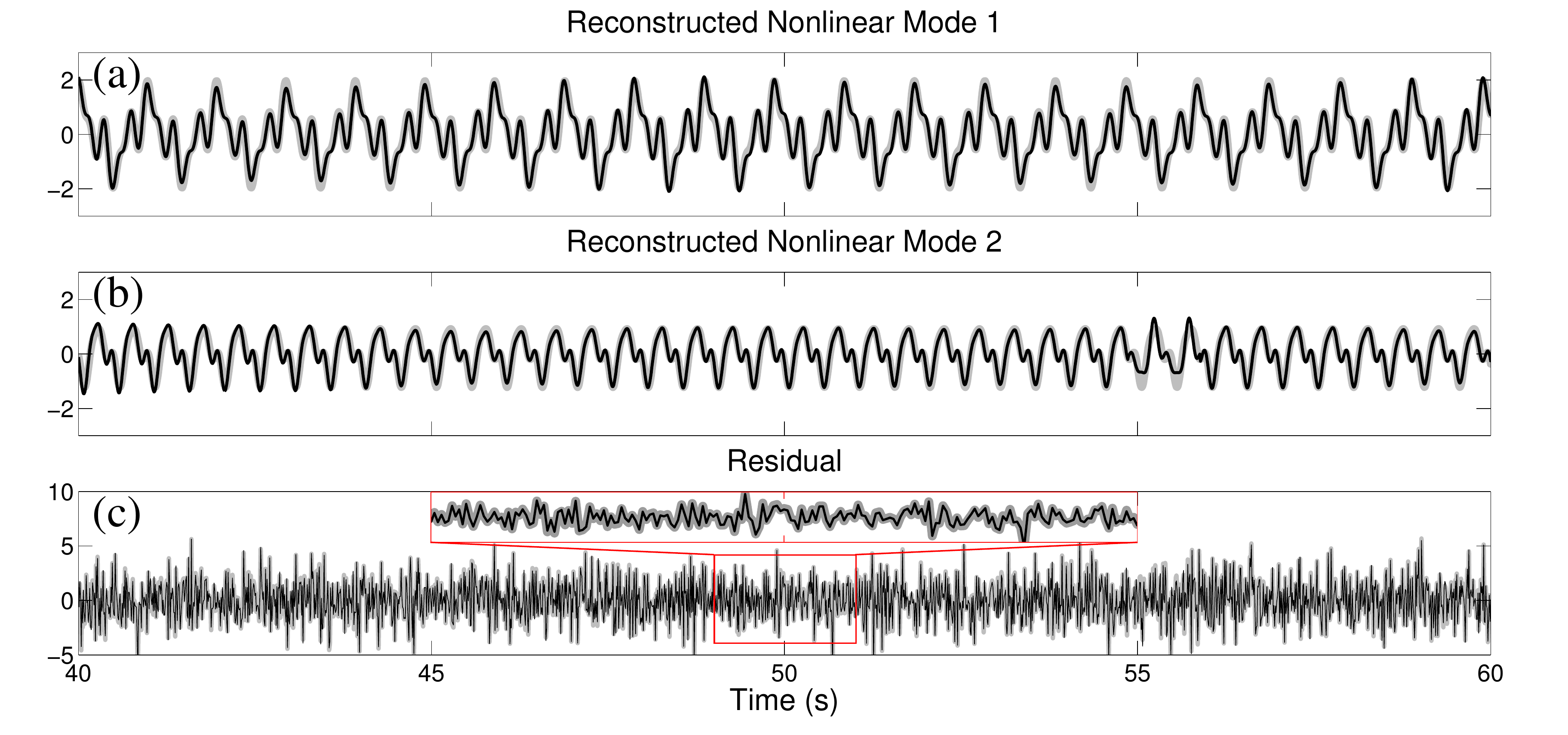}
\caption{The result of applying NMD to the signal shown in Fig.\ \ref{fig:ex2tfr}. In (a) and (b) black lines indicate the two reconstructed NMs and the gray background lines show the true NMs for comparison. (c) Similarly, the black and gray lines show the residual returned by NMD and true background noise, respectively.}
\label{fig:ex2nmd}
\end{figure}

The final results of NMD are shown in Fig.\ \ref{fig:ex2nmd}. Both NMs are reconstructed with great accuracy, which is like a miracle given such strong noise corruption of the original signal (see Fig.\ \ref{fig:ex2tfr}); even the residual noise is recovered almost exactly. Such performance is unachievable with the other existing methods, e.g.\ (E)EMD in the present case produces 13 components, and none of them reliably represent even a single harmonic (not shown). Note, that NMD can produce even more accurate results if the resolution parameter was adjusted from the very beginning, i.e.\ for the first harmonics (and not only the higher ones). However, as mentioned in Sec.\ \ref{sec:prespar}, at present there does not seem to be a good and universal way of doing this.

\section{Real life examples}\label{sec:realex}

After demonstrating its success in the analysis of simulated signals, we now apply NMD to real data.

\subsection{Decomposing human blood flow signals}\label{sec:BF}

The decomposition of skin blood flow signals, measured non-invasively by laser-Doppler flowmetry (LDF) \cite{Nilsson:80}, is believed to be a very tough task, with no method at present being able to do it well. This was demonstrated in \cite{Hozic:00} for the example of Karhunen-Lo\`{e}ve decomposition, while (E)EMD usually also fails. Thus, blood flow signals contain large amount of ``physiological'' noise, as well as having components with amplitude/frequency modulation whose strength and speed change with time. Nonetheless, as will be seen, NMD can often tackle even such complicated cases.

Blood flow signals exhibit oscillatory activity at multiple frequencies, and the WT has been found especially useful in studies of their structure \cite{Stefanovska:99a}. Six oscillations have been identified in blood flow and attributed to different physiological mechanisms \cite{Shiogai:10,Stefanovska:99a,Stefanovska:99b,Kvandal:06}, with characteristic frequency ranges of (approximately): 0.6-2 Hz (I), 0.15-0.6 Hz (II), 0.04-0.15 Hz (III), 0.02-0.04 Hz (IV), 0.01-0.02 Hz (V) and 0.005-0.01 Hz (VI). Range I corresponds to cardiac oscillations, which originate from heart activity, rhythmically pushing blood through the circulatory system, thus directly propagating into blood pressure and flow. Range II corresponds to respiratory oscillations, which are the consequence of the mechanical influence of respiration on the cardiac output and, to a smaller extent, the respiratory modulation of the heart rate \cite{Saul:91,Japundzic:90}. The III range oscillations by some authors are attributed to the sympathetic nerve activity, being regarded as a result of time-delays in the baroreflex feedback loop \cite{Malpas:02,Julien:06}, while the others relate them to myogenic activity of the smooth muscle cells \cite{Shiogai:10,Stefanovska:99a,Stefanovska:99b}. Finally, the oscillations in the IV, V and VI ranges were attributed to the neurogenic, NO-dependent endothelial and NO-independent endothelial activity \cite{Shiogai:10,Stefanovska:99a,Stefanovska:99b,Kvandal:06}, respectively.

Armed with this knowledge, we now apply NMD. To better utilize the prior information, we apply the procedure to each of the mentioned physiological frequency ranges individually, starting from the first one. Thus, for a given range we extract the dominant component and test it against noise; if it passes the test, we then extract its harmonics, reconstruct the NM and subtract it from the signal; we then repeat the procedure for the same range until the next extracted component does not pass the test against noise.

\begin{figure}[t]
\includegraphics[width=1.0\linewidth]{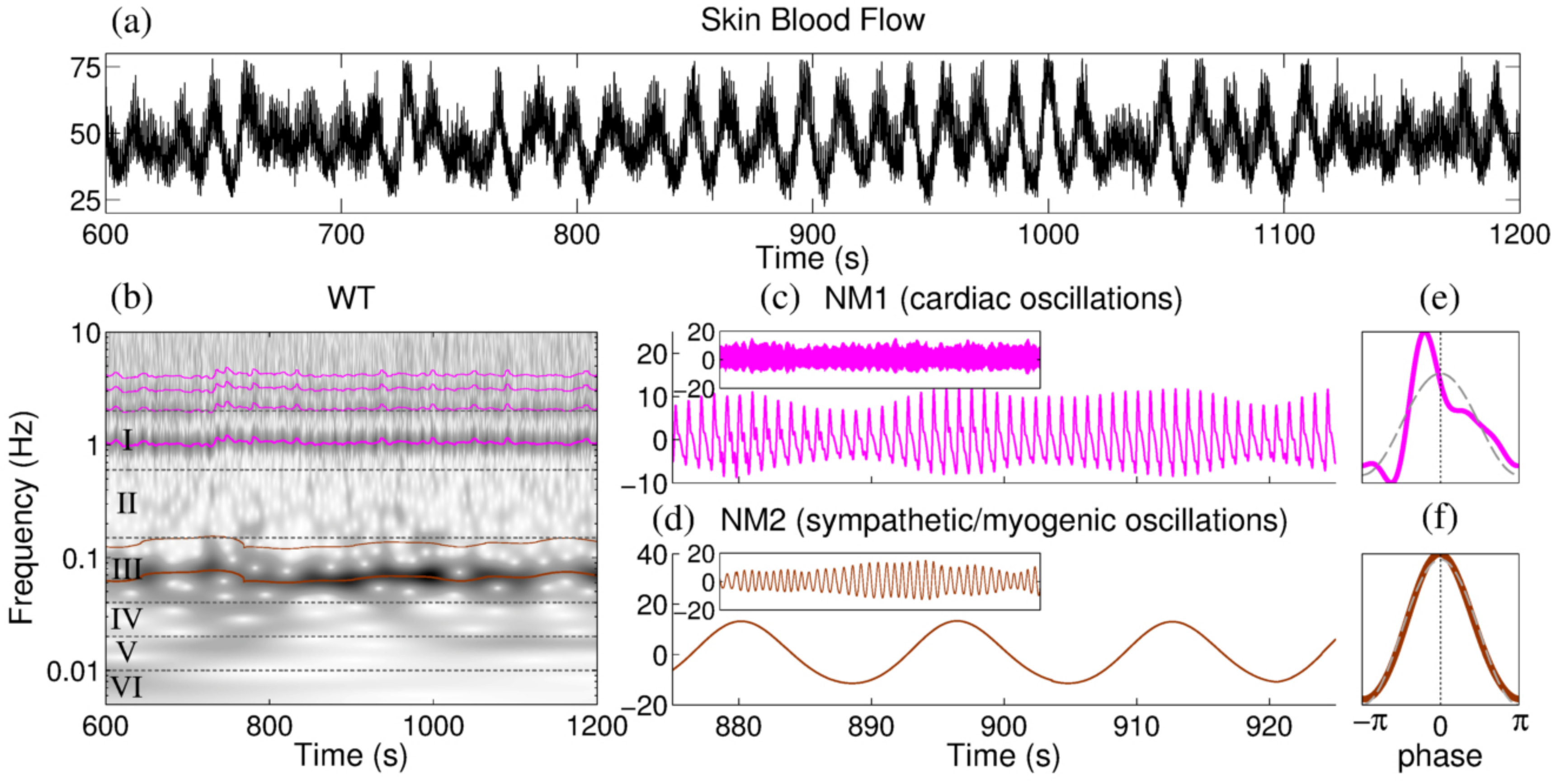}
\caption{(a) An example of the blood flow signal, measured by LDF with the probe over the right wrist \emph{caput ulna} (for more details, see \cite{Stefanovska:99b,Shiogai:10}); the signal is sampled at 40 Hz for 1800 s, and the panel shows its central 600\,s part. (b) The WT of the signal shown in (a). Gray-dotted lines partition the frequency axis into the regions within which the physiologically meaningful oscillations are located (according to \cite{Shiogai:10,Stefanovska:99a,Stefanovska:99b,Kvandal:06}, see the text). Bold colored lines show those extracted components which pass the surrogate test against noise, with thinner lines of the same color showing their higher harmonics. (c),(d) The reconstructed Nonlinear Modes, with the main graph showing them during 50 s, and small insets -- during 600 s (as in (a)); the colors of the lines correspond to those of the curves in (b). (e),(f) The waveforms of the oscillations shown in (c) and (d), respectively; the gray dashed line shows a pure sinusoidal waveform of amplitude equal to that of the first harmonic, and is provided as a guide to the eye. The cardiac waveform (e) has four harmonics $h=1,2,3,4$ with $a_h=[1,0.52,0.37,0.16]$ and $\varphi_h/\pi=[0,0.4,0.66,-0.97]$ (in the notation of (\ref{nm})); the waveform in (f) is characterized by $h=1,2$, $a_h=[1,0.08]$ and $\varphi_h/\pi=[0,0.07]$.}
\label{fig:BF112}
\end{figure}

\begin{figure}[t]
\includegraphics[width=1.0\linewidth]{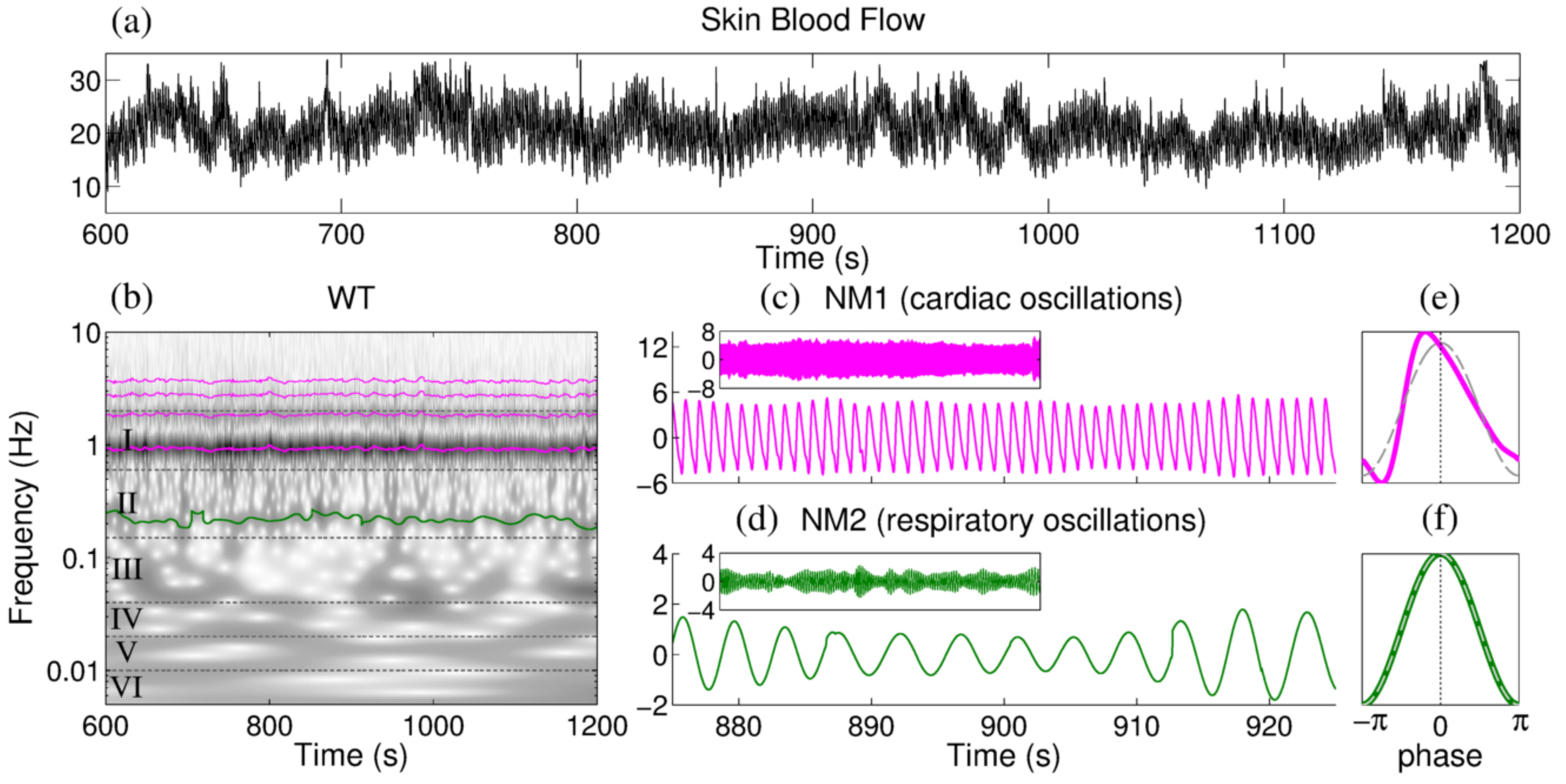}
\caption{Same as in Fig.\ \ref{fig:BF112}, but for the blood flow measured from a different subject. The cardiac waveform (e) has four harmonics $h=1,2,3,4$ with $a_h=[1,0.32,0.15,0.04]$ and $\varphi_h/\pi=[0,0.41,0.95,-0.44]$ (in the notation of (\ref{nm})); the respiratory oscillations have only one harmonic, so that the waveform in (f) is a pure sinusoid.}
\label{fig:BF108}
\end{figure}

The results of the procedure are shown in Figs.\ \ref{fig:BF112} and \ref{fig:BF108} for the examples of two blood flows. Clearly, NMD is able to decompose these signals into a physically meaningful oscillations with complex waveforms (and it also returns their amplitudes, phases and frequencies). In both cases, we were able to extract the cardiac component (around 1 Hz), while activity in ranges IV, V and VI did not pass the test against noise. However, in the example of Fig.\ \ref{fig:BF112} there are strong sympathetic/myogenic oscillations (around 0.1 Hz), which were extracted, while we see no activity in the frequency range II (respiratory) (Fig.\ \ref{fig:BF112}(b)). In contrast, for the example of Fig.\ \ref{fig:BF112} the respiratory oscillations are present and the sympathetic/myogenic are buried under noise. The waveforms of the cardiac oscillations are also different in two cases.

In general, the properties and presence of the oscillations in blood flow varies from subject to subject, being influenced by many factors, such as the state of the microvasculature (which might in turn be influenced by age and gender), properties of the skin etc. As was demonstrated, NMD can be very useful for the study and classification of these effects, e.g.\ it can be used to investigate the relationship between the cardiac waveform (as well as other blood flow properties) and the health conditions.

\begin{remark}
As discussed in Sec.\ \ref{sec:stopcrit}, even if the component extracted from a particular frequency range does not pass the surrogate test against noise (as e.g.\ for ranges IV, V and VI in the above examples), this does not necessarily mean that there are no physiologically meaningful activity there. Thus, the underlying oscillations might be simply very small so that they are easily masked by noise, which is often the case for the respiratory oscillations. The other possibility is that the resolution parameter used is not appropriate to represent reliably the component of interest. This is often the case for sympathetic/myogenic oscillations, which might change their amplitude/frequency very rapidly at certain times. In fact, the best choice in such situations would be probably some time-varying $f_0(t)$, but its form is generally very hard to choose.
\end{remark}

It should be noted, that NMD can also serve as an initial preprocessing which needs to be performed prior to applying any of the numerous existing methods for studying monocomponent signals or their sets. Thus, having first decomposed the original signal into its constituent NMs, one can then investigate the latter using one or some of the huge diversity of available techniques. For example, the structure of the interactions between different modes can be recovered by applying Bayesian inference to the extracted phases \cite{Stankovski:12,Duggento:12}, as was done in \cite{Iatsenko:13a} to reconstruct the cardiorespiratory interaction; in this application, the high accuracy of the phases returned by NMD is especially advantageous.

Another important application is the classification of the oscillations present in the signal, which is of particular interest because it might yield valuable insights into the nature and properties of the underlying phenomena. Following the recent introduction of chronotaxic systems \cite{Suprunenko:13,Suprunenko:14,Clemson:14}, it is clearly desirable to be able to determine whether the originating systems generating the different modes are chronotaxic or not. Roughly speaking (see \cite{Suprunenko:13,Suprunenko:14} for a more detailed definition), the system is chronotaxic if it is: (a) oscillatory (i.e.\ characterized by a limit cycle); and (b) its phase $\phi(t)$ does not just move freely along the cycle, as conventionally assumed, but is attracted to some $\phi_A(t)$, conferring an ability to resist external perturbations. This is exactly what one often observes in living systems, which are able to maintain their activity within physiological ranges even when strongly perturbed, so that chronotaxic behavior is expected to be abundant within the life sciences.

To establish which oscillations contained in the signal are chronotaxic, one needs to study them separately, i.e.\ the signal should be first decomposed into its NMs; this can conveniently and accurately be achieved with NMD. Subsequently, for each mode one can apply the approach suggested by Clemson et al.\ \cite{Clemson:14}. Basically, by applying different types of filters to the extracted phase one estimates the expected difference $\phi(t)-\phi_A(t)$ and uses detrended fluctuation analysis \cite{Peng:94,Shiogai:10} to analyze the associated correlations, which are expected to differ between chronotaxic and non-chronotaxic systems. We have applied this method to the NMs of both of the examples in Figs.\ \ref{fig:BF108} and \ref{fig:BF112}, but have not found clear evidence of chronotaxicity in any of the corresponding oscillations. However, since the method of Ref.\ \cite{Clemson:14} is based on a particular set of assumptions (e.g.\ that the dynamical perturbations take form of the Gaussian white noise), these oscillations could in principle still be chronotaxic though falling outside the model considered in \cite{Clemson:14}. 

\subsection{Removing cardiac artifacts from a human EEG signal}\label{sec:EEG}

Nonlinear Mode Decomposition can also be used to filter the signal from an extraneous oscillatory activity, provided that there is an associated signal from which the phase and frequency of the latter can accurately be extracted. Thus, using the NMD harmonic extraction procedure, the fundamental harmonic of the extraneous mode can be extracted as $h=1$ harmonic of the reference component based on its (reference) phase and frequency. In this way the resolution parameter is adjusted from the start, allowing one to represent this fundamental oscillation well even if it is strongly corrupted by noise or other influences. We will illustrate this alternative use of NMD through the example of removing cardiac artifacts from the human electroencephalogram (EEG) recording, using the electrocardiogram (ECG) as the reference signal.

The EEG often contains artifacts related to heart activity, which arise due to blood flow pulsations underneath the probes (the so-called ballistocardiogram (BCG) artifacts \cite{Srivastava:05}) and possibly also due to direct pick up of heart electrical activity. The BCG artifacts are extremely prominent in the EEG measured in a magnetic field (e.g.\ simultaneously with a magnetic resonance imaging scan), in which case they are usually filtered out by independent component analysis (ICA) \cite{Comon:94,Hyvarinen:00}. However, ICA requires as many simultaneous EEG measurements as possible, with cardiac artifacts being prominent and taking a similar form in most of them. Therefore, it will obviously fail to remove the artifacts if only one EEG and one ECG signal are available, as the form taken by cardiac activity in the EEG is completely different from its form in the ECG signal (see below), contradicting the assumption of linear mixing on which ICA is based. In fact, we have found that ICA fails even given four EEGs, which is probably because the artifacts in them are relatively small, being hard to distinguish (although still capable of affecting some time-frequency measures); and, additionally, because the form of these artifacts might be different in different EEGs (perhaps dependent on probe position). The (E)EMD method, too, fails to provide meaningful results in the present case.

\begin{figure}[t!]
\includegraphics[width=1.0\linewidth]{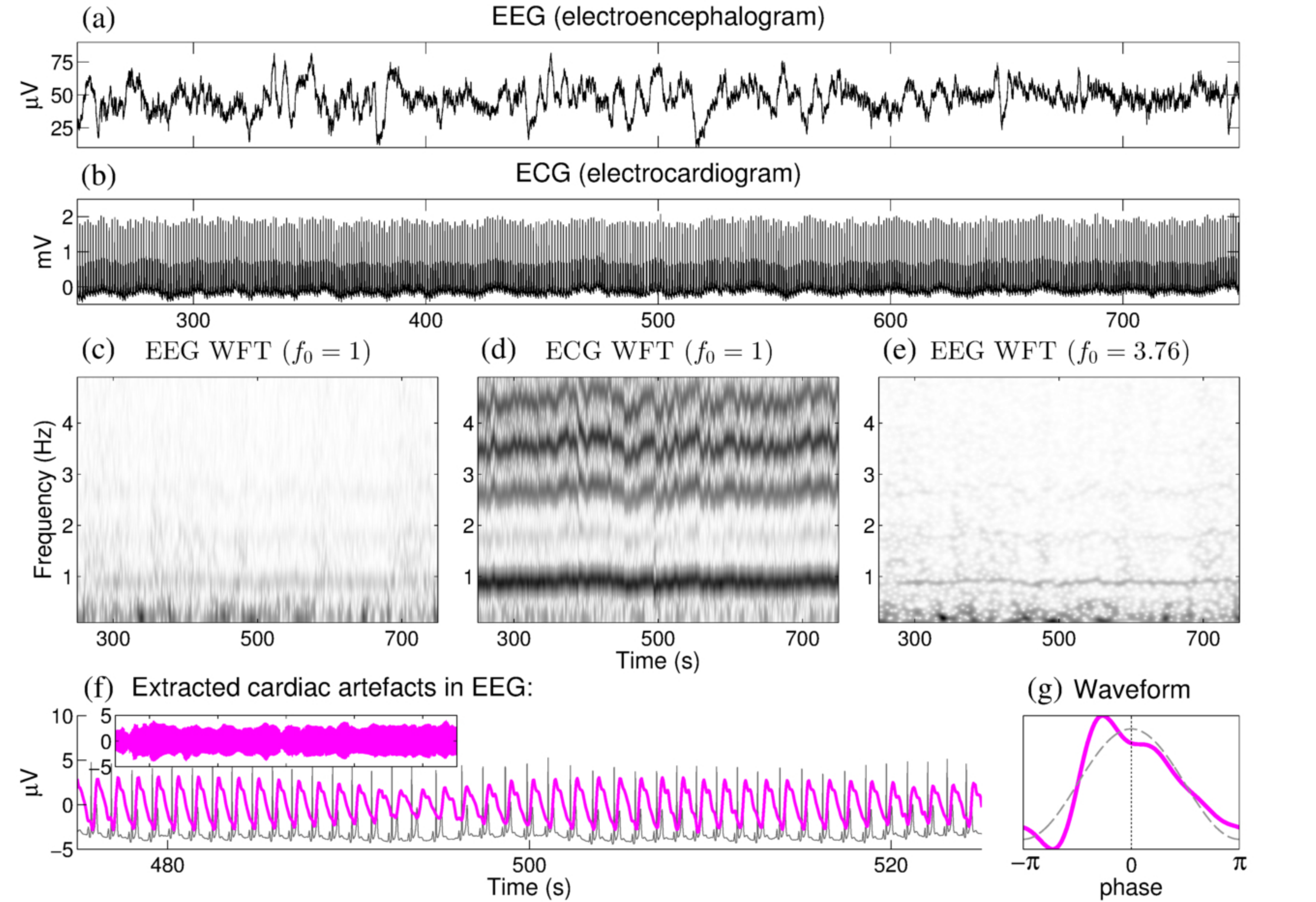}\\
\caption{(a) The EEG signal (measured using BIS electrode placed on the forehead, as described in \cite{BRACCIA}). (b) The ECG signal (3-lead, with electrodes on shoulders and the lowest left rib, see e.g.\ \cite{Lotric:00}). (c),(d) The WFTs of the EEG and ECG signals shown in (a) and (b), respectively, calculated using the default resolution parameter $f_0=1$. (e) The WFT of the EEG signal calculated using the adaptively adjusted resolution parameter $f_0=3.76$, due to which the cardiac component becomes much more visible than in (c). (f) The cardiac artifacts extracted from the EEG, with an inset showing them over the same 500 s as presented in (a,b); gray background lines show the ECG scaled to the dimensions of the plot. (g) The waveform of the cardiac artifacts, which has four harmonics with $a_h=[1,0.33,0.24,0.11]$ and $\varphi_h/\pi=[0,0.49,1,-0.59]$ (in the notation of (\ref{nm})); the gray-dashed line represents a pure sinusoid with the amplitude of the first harmonic, for comparison. The full signals are of duration 20 min.\ and are sampled at 80 Hz; the EEG was actually measured under anaesthesia, but the same artifacts arise under all conditions.}
\label{fig:EEGnmd}
\end{figure}

Simultaneously measured EEG and ECG time series for the same subject are presented in Fig.\ \ref{fig:EEGnmd}(a) and (b), with their WFTs for the default resolution parameter $f_0=1$ being shown in (c) and (d), respectively. Clearly, for $f_0=1$ the cardiac harmonics are not well distinguishable in the EEG's WFT (Fig.\ \ref{fig:EEGnmd}(c)), so that the extracted curve might not be very accurate. However, the cardiac phase $\phi^{(1)}(t)$ and frequency $\nu^{(1)}(t)$ can of course be extracted {\it directly} from the WFT of the ECG (Fig.\ \ref{fig:EEGnmd}(d)). They should be the same for the first cardiac harmonic in the EEG and in ECG, because both activities obviously have the same rhythmicity; however, the corresponding amplitudes might be different, and perhaps not even proportional to each other. This is because, depending on the measurement and the environment, the same activity might undergo various transformations that can change its amplitude dynamics and the corresponding waveform, but leave phase dynamics largely unaltered. For example, Nonlinear Modes $c(t)=A(t)[\cos\phi(t)+a\cos2\phi(t)]$ and $\exp[c(t)]$ will have the same fundamental phase and frequency, but different amplitude dynamics (i.e.\ the ratio between the amplitudes of the corresponding fundamental components will be time-varying) and different relationships between the harmonics.

Therefore, we use the cardiac phase $\phi^{(1)}(t)$ and frequency $\nu^{(1)}(t)$ extracted from the ECG as the parameters of the reference component, and extract the main cardiac component from the EEG as its $h=1$ harmonic. This is done in the same way as discussed in Sec.\ \ref{sec:hextract}, i.e.\ selecting the peaks nearest to the expected frequency $\nu^{(1)}(t)$. We test this harmonic for being true (Sec.\ \ref{sec:htest}) and adapt the resolution parameter for its best representation (Sec.\ \ref{sec:ihadapt}) in the usual way, except that now we use everywhere the phase consistency $\rho^{(1)}(0,1,0)$ (\ref{apc}) instead of the usual amplitude-phase consistency $\rho^{(1)}(1,1,0)$. This is because of the possible discrepancy between the amplitude dynamics of the components mentioned above (note that the threshold (\ref{rhomin}) also becomes $\rho_{\min}=0.5$). If the extracted first harmonic is regarded as true, we use its parameters to extract the higher cardiac harmonics from EEG (now using unmodified procedures) and reconstruct the corresponding NM.

\begin{remark}
Note that there might in general be a time delay between related activities in different signals, e.g.\ ECG and cardiac artifacts in EEG. However, unless it exceeds the characteristic time of amplitude/frequency variations or the minimal surrogate time shift (see Sec.\ \ref{sec:htest}), it should not influence the results significantly. Nevertheless, to be fully rigorous one can additionally adjust the timings of both components by maximizing frequency consistency $\rho^{(1)}(0,0,1)$ (\ref{apc}) between them.
\end{remark}

The EEG's WFT adapted (by maximizing phase consistency $\rho^{(1)}(0,1,0)$) for representation of the cardiac component is shown in Fig.\ \ref{fig:EEGnmd}(e). Clearly, the corresponding ridge curves became much more visible than in the default WFT presented in Fig.\ \ref{fig:EEGnmd}(c). The cardiac artifacts extracted from the EEG are shown in Fig.\ \ref{fig:EEGnmd}(f). Their waveform, presented in Fig.\ \ref{fig:EEGnmd}(g), very much resembles the shape of the cardiac waves in the blood flow (cf.\ Fig.\ \ref{fig:BF108}), but not that of the ECG. This is an indication of the BCG mechanism by which these cardiac artifacts are generated, which is also supported by the fact that their strength (and even their shape) might be different in EEGs from different probes for the same subject. Note that, depending on the particular EEG measurement, the artifacts might be inverted.

There are many other possible applications of the discussed NMD-filtering. For example, given the reference ECG and respiration signals, one can use the same approach to extract the heart and respiratory components from the blood flow more accurately (especially for respiratory activity) than is possible using the straightforward decomposition as in Sec.\ \ref{sec:BF}. Note that if in the given signal there is no oscillation related to the reference phase and frequency, then the corresponding first harmonic extracted from it will be regarded as a false harmonic of the reference component. Thus, we have not found any respiratory oscillations in the EEG, implying that the measurement process is almost unaffected by breathing.

In general, the situation when one signal contains components related to other signals is ubiquitous in real life, so that the NMD-based filtering is expected to be very useful in many contexts. A great advantage of this approach is that it does not require related oscillations in different signals to be of the same form (as is the case e.g.\ for ICA), but only to have the same phase dynamics.

\section{Conclusions}\label{sec:conclusions}

We have introduced a new decomposition method -- Nonlinear Mode Decomposition (NMD), which is based on the time-frequency analysis techniques \cite{Iatsenko:13tfr1,Iatsenko:13tfr2,Iatsenko:13ridge}, surrogate data tests \cite{Schreiber:00,Theiler:92,Andrzejak:03,Quiroga:02}, and the idea of harmonic identification \cite{Sheppard:11}. NMD consists of a number of ingredients, as described in Sec.\ \ref{sec:nmd}, each of which is a useful technique in its own right. For example, one can use the approach discussed in Sec.\ \ref{sec:stopcrit} for testing the signal against noise, which proves to be much more noise-robust than the conventional tests based on the correlation dimension or other nonlinearity measures \cite{Schreiber:97,Schreiber:00,Theiler:92}.

We have successfully applied NMD to both simulated and real data, and demonstrated its advantages over the existing approaches. Unlike many previous methods, NMD recovers the full underlying oscillations of any waveform. It is also extremely noise-robust and in a sense super-adaptive, meaning that most of its settings are automatically adapted to the properties of the particular signal under investigation. Finally, in contrast to the other methods, NMD returns only physically meaningful modes, stopping the decomposition when the residual is just noise.

The area of applicability of NMD is extremely wide. For example, it would now be reasonable to reconsider all those cases where (E)EMD has been applied to date (e.g.\ see references in \cite{Huang:08,Wu:09}). Furthermore, as demonstrated above, the exceptional noise-robustness of NMD and its other advantages allow one to apply it even in cases where all the other methods fail completely. Thus, it can be applied to almost all multicomponent signals of the kind that are ubiquitous in the life sciences, geophysics, astrophysics, econometrics etc. We therefore expect that, with time, NMD will become a new standard in the field.

The latest MATLAB codes needed for running NMD are freely available \cite{FreeCodes}, together with detailed instructions and examples, in both text and video formats.

\section*{Acknowledgement}

We are grateful to Philip Clemson for valuable discussions. The work was supported by the Engineering and Physical Sciences Research Council (UK) [grant number EP/100999X1].


\bibliographystyle{elsarticle-num}    %
\bibliography{NMDbib}

\end{document}